\documentclass[11pt]{article}
\usepackage{times,amsmath,amsfonts,amsthm,amssymb,eucal}
\usepackage{algorithm}
\usepackage{algorithmic}
\usepackage{wrapfig}
\usepackage{color}

\usepackage{natbib,graphicx}

\newcommand\BibTeX{{\rmfamily B\kern-.05em \textsc{i\kern-.025em b}\kern-.08em
T\kern-.1667em\lower.7ex\hbox{E}\kern-.125emX}}

\newtheorem{theorem}{Theorem}[section]
\newtheorem{lemma}{Lemma}[section]
\newtheorem{remark}{Remark}[section]
\newtheorem{proposition}{Proposition}[section]

\begin{document}

\begin{center}

\Large{\bf Adaptation in some linear inverse problems}
%\title{Adaptation in a class of linear inverse problems}

\vskip.1in
\large{Iain M. Johnstone$^\dag$ and Debashis Paul$^\ddag$}

\vskip.1in
\large{\textit{$\dag$ Stanford University, $\ddag$ University of California, Davis}}

\vskip.15in \large{\bf This paper is dedicated to Laurent Cavalier}

\end{center}

%\address{%
%\affilnum{a}Department of Statistics, Stanford University\\
%\affilnum{b}Department of Statistics, University of California, Davis
%}

%\corremail{imj@stanford.edu}

%\received{00 Month 2014} \accepted{00 Month 2014}

\begin{abstract}
We consider the linear inverse problem of estimating an unknown signal $f$ from
noisy measurements on $Kf$ where the linear operator $K$ admits a
wavelet-vaguelette decomposition (WVD).
%This refers to a setup which can be described by the following sequence model:
%$$%y_{jk} = \lambda_j \theta_{jk} + \epsilon z_{jk},\qquad k=0,1,\ldots
%2^j-1, \quad  j=1,2,\ldots
%$$%where $\theta = (\theta_{jk})$ is the vector of wavelet coefficients
%of the unknown function $f$. The weights $\lambda_j$
%satisfy $\lambda_j \sim 2^{-\alpha j}$ for some $\alpha > 0$. And the noise
%$z_{jk}$ are assumed to be Gaussian with 0 mean, and their covariance matrix
%has bounded condition number. We proceed to estimate $\theta$ (equivalently
%$f$)
We formulate the problem in the Gaussian sequence model and propose
estimation based on complexity penalized regression
%scheme for the vaguelette coefficients of the signal
on a level-by-level basis.
We adopt squared error loss and show that the estimator achieves exact
rate-adaptive optimality
% is used as performance measure for estimation of $f$, then the
% proposed estimator achieves the optimal rate of convergence in the small noise
% limit
as $f$ varies over a wide range of Besov function classes.
\end{abstract}


%\maketitle

\section{Introduction}\label{sec:introduction}

%\textcolor{red}{[TO DO: Change this section by bringing in the main features of
%the problem and giving a summary of the results as in Iain's draft. Also, make
%it clear that numerical simulations and realistic applications are not being
%considered owing to space constraints. Follow the structure suggested by IMJ.]}

%\vskip.15in

This paper studies the recovery of an unknown function $f$ based on noisy
measurements on $g = Kf$ where $K$ is a linear operator belonging to a class of
homeogeneous, ill-posed operators. To set the stage, we recall the
\textit{direct} estimation setting, $Y_\epsilon (dt) = f(t) dt + \epsilon
W(dt)$, where $W(t)$ is the standard Browian motion. Here it is now well
understood that expansion in wavelet bases is useful for the estimation of
spatially inhomogeneous functions $f$. Spatial inhomogeneity is formulated by
supposing that $f$ belongs to an appropriate Besov space. Wavelet shrinkage
estimators are shown to have adaptive minimaxity properties over a wide range
of Besov space function classes. The key property of adaptivity means that the
wavelet estimator attains the optimal rate of convergence for each Besov class
even though the estimator is specified without knowledge of the parameters of
that Besov class.

The goal of this paper is to exhibit the first wavelet-type estimator with exact
rate-adaptive optimality in a class of ill-posed linear inverse problems where
the observed data can be described through the model
\begin{equation}\label{eq:white_noise}
Y_\epsilon(dt) = (Kf)(t) dt +  \epsilon W(dt), \qquad t \in [0,1]
\end{equation}
and $K$ is a linear operator acting on ${\cal D}(K) \subset L_2([0,1])$. The
inverse problems we consider are those in which $K$ possesses a
wavelet-vaguelette decomposition (WVD), to be recalled below.
Homogeneous operators, which satisfy $(Kf)(at) =
a^{-\beta}(Kf)(t)$ for all $t$, $a>0$ and some $\beta \in \mathbb{R}$,
provide a class of examples, such as $r$-fold integration for arbitrary positive integer $r$,
fractional integration and convolution with a suitably regular convolution
kernel. A two-dimensional example
%of such operators when the domain is two-dimensional
is the Radon transform, seen for example in positron emission
tomography \citep{Kolaczyk96}.
%considered a wavelet-vaguelette decomposition of the Radon
%transform.}

A detailed discussion of the conditions on the operators and the
function spaces involved can be found in \cite{Donoho95}, while
\cite{JohnstoneKPR04} treated in detail the case where $K$ is a convolution
operator, i.e., $Kf = K*f$ for a kernel $K \in L_1([0,1])$.  These cases are
characterized by an ill-posedness index $\beta \geq 0$ and, as recalled below,
can be recast, using the WVD, into the form of a Gaussian multi-resolution
sequence model
\begin{displaymath}
y_{jk} = \theta_{jk} + \epsilon 2^{\beta j} z_{jk}.
\end{displaymath}
Here the ill-posedness index $\beta$ appears as a noise inflation factor.

%\subsection{Minimax rates over Besov spaces}
The function classes we consider are represented by Besov norm balls
$\Theta_{p,q}^\alpha (C)$ indexed by smoothness $\alpha$ and radius $C$. Here
$p$ is the $L_p$ integration parameter, with $p < 2$ corresponding to cases
modeling spatial inhomogeneity. The minimax mean squared error for such a class
$\Theta$ is
\begin{displaymath}
  R_N (\Theta, \epsilon)
   = \inf_{\hat \theta} \sup_{\theta \in \Theta}
        E_\theta \| \hat \theta - \theta \|^2.
\end{displaymath}

The rate of convergence will be defined by a ``rate control function'' $R(C,
\epsilon; \boldsymbol{\gamma})$ depending on a vector parameter
$\boldsymbol{\gamma} = (\alpha, p, q , \beta)$ ranging over a set $\Gamma$. The
vector $\gamma$ encodes the parameters of the Besov ball and the ill-posedness
index $\beta$. The form of $R(C, \epsilon; \gamma)$ depends on one of the three
zones comprising $\Gamma$ to be described later; for example, for
$\boldsymbol{\gamma}$ in the ``dense'' zone $\Gamma_d$ we have
\begin{displaymath}
  R(C, \epsilon; \boldsymbol{\gamma}) = C^{2(1-r)} \epsilon^{2r},
    \qquad r = 2 \alpha/(2 \alpha + 2 \beta + 1).
\end{displaymath}
The rate therefore depends on both the smoothness $\alpha$ and the
ill-posedness index $\beta$, with increasing $\beta$ leading to slower rates of
convergence.

The main result of this paper---stated more formally in Theorem
\ref{thm:minimax_risk} below---is the construction of a penalized least squares
estimator $\hat \theta_P$ and the demonstration that it satisfies, for
$\boldsymbol{\gamma} \in \Gamma$ and all $\epsilon$ sufficiently small,
\begin{equation*}
  \begin{split}
  c_0 R(C, \epsilon; \boldsymbol{\gamma})
    & \leq R_N( \Theta_{p,q}^\alpha (C), \epsilon)  \\
    & \leq \sup_{\Theta_{p,q}^\alpha (C)} E \| \hat \theta_P - \theta
    \|^2
       \leq c_1 R(C, \epsilon; \boldsymbol{\gamma}) + c_2 \epsilon^2 \log
       \epsilon^{-2}.
  \end{split}
\end{equation*}
The term $\epsilon^2 \log \epsilon^{-2}$ is of smaller order than $R(C,
\epsilon; \boldsymbol{\gamma})$. The constants $c_1$ and $c_2$ depend on
$\gamma$ but the chief conclusion is that $\hat \theta_P$ achieves the minimax
rate of convergence for each $\boldsymbol{\gamma} \in \Gamma$, and does so
without knowledge of $\boldsymbol{\gamma}$.

%Suppose we are interested in recovering an unknown function $f$ that is a
%member of $L^2(\mathbb{R}^d)$ when we are able to observe measurements on $g =
%Kf$ where $K$ is a known linear operator. Examples of such operators include
%the Radon transform, convolutions and fractional integration.
There is an extensive literature on linear inverse problems
in statistics. One may refer to
\cite{AbramovichS98}, \cite{BissantzHMR07}, \cite{Cai02}, \cite{CavalierG06},
\cite{CavalierGLT04},  \cite{CavalierGPT02}, \cite{CavalierR07},
\cite{CavalierT02}, \cite{Rochet13},
\cite{Donoho95}, %\citeN{FanK02},
\cite{Johnstone99}, \cite{JohnstoneKPR04}, \cite{KalifaM03}, \cite{Kolaczyk96},
\cite{LoubesL08} and \cite{PenskyV97}, among others, for some recent
advances in this field. Specifically, \cite{Donoho95} proposed solving the
linear inverse problems described above through the WVD framework and obtained
lower bound on the rate of convergence in the ``dense'' regime. He also
proposed an estimator that can attain the optimal rate of convergence under the
$L_2$ loss, with the knowledge of the hyperparameters of the Besov function
class. \cite{CavalierR07} considered an estimator based on hard thresholding of
the empirical Fourier coefficients and derived upper bounds for the rate of
convergence under $L_2$ loss, which are within a factor of $\log \epsilon^{-2}$
of the optimal rate in all three (``dense'', ``sparse'' and ``critical'')
regimes. \cite{Cavalier08} and \cite{LoubesR09} gave nice surveys of the
various approaches to statistical inverse problems and summarized results on
the rates of convergence of the estimators.

We now review aspects of the WVD and its relation to our sequence model.
Given a wavelet basis with mother wavelet $\psi$ and scaling function $\phi$,
satisfying appropriate regularity conditions, there are \textit{biorthogonal
systems of vaguelettes} ${\cal U}$ and ${\cal V}$ and a sequence of
pseudo-singular values $\kappa_j$ (depending on the scale index $j$ but not on
the spatial index $k$) such that, formally,
\begin{equation}\label{eq:vaguelettes}
K\psi_{jk} = \kappa_j v_{jk}, \qquad Ku_{jk} = \kappa_j \psi_{jk},
\qquad\mbox{and hence}~~ \langle u_{jk}, v_{j'k'}\rangle = \delta_{j-j'}
\delta_{k-k'}
\end{equation}
where $\delta_j$ denotes the Kronecker's delta function. Supposing that we have
a representation of the function $f$ in the inhomogeneous wavelet basis as
$$
f(t) = \sum_{k} \langle f,\phi_{j_0 k}\rangle \phi_{j_0 k}(t) +
\sum_{j=j_0}^\infty \sum_k \langle f,\psi_{jk}\rangle \psi_{jk}(t)
$$
we can use (\ref{eq:vaguelettes}) to write $\theta_{jk} := \langle
f,\psi_{jk}\rangle = \kappa_j^{-1} \langle Kf, u_{jk}\rangle$. The coefficients
$\langle f,\phi_{j_0 k}\rangle$ can be obtained as a linear combination of the
coefficients $\theta_{j_0 k}$ since the function $\phi_{j_0 0}$ can be
expressed as a linear combination of the functions $\{\psi_{j_0 k}\}$.

The \textit{frame property} of the WVD system \citep{Donoho95} states that
there exist constants $0 < \Xi_0 < \Xi_1 < \infty$ so that
\begin{equation}\label{eq:frame_u_v}
\Xi_0 \parallel (\alpha_{jk}) \parallel_2^2 ~\leq~  \parallel \sum_{j,k}
\alpha_{jk} u_{jk} \parallel_2^2 ~\leq~ \Xi_1 \parallel (\alpha_{jk})
\parallel_2^2
~\mbox{and}~~ \Xi_0 \parallel (\alpha_{jk}) \parallel_2^2 ~\leq~  \parallel
\sum_{j,k} \alpha_{jk} v_{jk} \parallel_2^2 ~\leq~ \Xi_1
\parallel (\alpha_{jk}) \parallel_2^2
\end{equation}
for any sequence $(\alpha_{jk}) \in \ell_2$. In other words, if (by an abuse of
notation) we denote by ${\cal U}$ and ${\cal V}$ the operators corresponding to
the vaguelette transform on appropriate domains, then the Gram operators
satisfy
\begin{equation}\label{eq:Gram}
\Xi_0 I ~\leq~ {\cal U}^*{\cal U} ~\leq~ \Xi_1 I, \qquad \qquad\mbox{and}\qquad
\qquad \Xi_0 I ~\leq~ {\cal V}^*{\cal V} ~\leq~ \Xi_1 I.
\end{equation}

Therefore one may perform a vaguelette transform of the data in model
(\ref{eq:white_noise}) in the system ${\cal U}$ and get empirical wavelet
coefficients $y_{jk} = \kappa_j^{-1}\langle Y,u_{jk}\rangle$.
We assume that $\sigma$ is known and set $\epsilon =
\sigma n^{-1/2}$.
% %\subsection{WVD setting}\label{subsec:WVD}
% We consider the white noise model (\ref{eq:white_noise}) and set $\epsilon =
% \sigma n^{-1/2}$ assuming that $\sigma$ is known. We apply a vaguelette
% transform to the data $Y$ to get the empirical vaguelette coefficients $y_{jk}
% = \kappa_j^{-1} \langle Y,u_{jk}\rangle$.
Then, we can express the model
(\ref{eq:white_noise}) in the transformed system as
\begin{equation}\label{eq:WVD_transform}
y_{jk} = \theta_{jk} + \epsilon_j z_{jk}, \qquad k=1,\ldots,2^j, \quad j\geq
j_0
\end{equation}
where $z_{jk} = \int_0^1 u_{jk}(t) dW(t)$ and $\epsilon_j = \epsilon \sigma_j$
with $\sigma_j = \kappa_j^{-1}$. Therefore the variance of the noise $z_{jk}$
at the dyadic level $j$ is $\epsilon_j^2$. Henceforth, we define $n_j = 2^j$.
Let $\Sigma_j$ denote the covariance matrix of $z_j = (z_{jk} :
k=1,\ldots,2^j)$, with $\parallel \Sigma_j \parallel ~= \xi_j$. Observe that
$\xi_j \leq \Xi_1$. In many cases, for example when $K$ is a convolution
operator, $\xi_j$ is numerically computable.

We assume that the pseudo-singular values
$\kappa_j$ of the operator $K$, or equivalently, their inverses $\sigma_j = \kappa_j^{-1}$, satisfy
\begin{equation}\label{eq:bound_singular}
B_0 \leq \sigma_j 2^{-\beta j} \leq B_1 \qquad \mbox{for some} ~~ \beta > 0
\end{equation}
and for constants $0 < B_0 \leq B$. To keep the exposition simple and avoid
cumbersome expressions, throughout we assume that $B_0 = B_1 = 1$, i.e., $\sigma_j =
2^{\beta j}$, since the discrepancy between $\sigma_j 2^{-\beta j}$ and 1 can
be absorbed in the covariance matrix $\Sigma_j$ and this simplification can
only change the extreme singular values of $\Sigma_j$ by a constant multiple.

Then, one can estimate $f$ by obtaining estimates of the coefficients $\theta_{jk}$ derived
by regularizing $y_{jk}$. Indeed, the central proposal of \cite{Donoho95}
%showed that when $f$ is known to belong to a Besov space,
%(with certain restrictions on the parameters)
is to estimate $f$ by coordinate-wise hard thresholding of the coefficients
$y_{jk}$, assuming that the function $f$ belongs to a certain Besov function
class. This estimator has asymptotically optimal rate of convergence over the
Besov function class in the minimax sense as $\epsilon \rightarrow 0$. However,
this estimator requires the knowledge of the hyperparameters defining the Besov
class to which $f$ belongs and therefore it is not adaptive. The estimator of
\cite{JohnstoneKPR04} has rate of convergence within a factor of (a power of)
$\log \epsilon^{-2}$ of the minimax rate, even though it does not require the
knowledge of the hyperparameters. \cite{Cai02} obtained similar results for a
level-wise James-Stein estimator of the vaguelette coefficients.

%Our aim here is to propose and analyze an estimator that is rate adaptive over
%a wide variety of Besov function classes.

%In this paper, we adopt a level-wise penalized least squares approach to solve
%this problem.
%Before going into the details of the method it is important to
%point to a crucial fact about the WVD system, that is its \textit{frame
%property}. This is a restatement of Theorem 2 of \cite{Donoho95}.

\textbf{Outline of the paper.} \
The frame property (\ref{eq:frame_u_v}) obviously
holds for each individual scale (i.e., for each $j$) and therefore we shall
first provide a general monoscale estimation procedure in a gaussian linear
model setup known as {\it additive, weakly correlated noise}. This is the topic
of Section \ref{sec:correlated_noise}. In Section \ref{sec:besov_minimax} we
deal specifically with the WVD paradigm and propose a multiscale estimation
procedure which uses a penalized estimator for each scale separately. In
Section \ref{sec:upper_bound} we produce upper bounds on the risk of our
estimator. In Section \ref{sec:lower_bound} we provide the matching lower
bounds that prove the rate-optimality of the proposed estimator. Owing to space
constraints, numerical simulations and realistic applications are not
considered in this paper.
In order to give a compact account of the derivations
in Sections \ref{sec:upper_bound} and \ref{sec:lower_bound}, we refer to well established material covered in
\cite{Johnstone2013}.
% Derivations of the main results are
% based on the material in Chapters 10, 11 and 12 of
% \cite{Johnstone2013}.
In addtion, some technical details are provided in the Supplementary
Material (SM).

%\vskip.1in\textcolor{red}{TO DO: In notations for Besov spaces, change $s$ to $q$.}

\section{Penalized estimation at one resolution level}\label{sec:correlated_noise}
%\section{Regression under WVD framework}\label{sec:correlated_noise}

We first consider estimating a sparse parameter vector in presence of additive,
weakly correlated Gaussian noise:
\begin{equation}\label{eq:general}
y_k = \theta_k + \epsilon z_k, \quad k=1,\ldots,n.
\end{equation}
Here, the noise vector $Z = (z_k)_1^n$ is distributed as $N(0,\Sigma)$. Let
$\xi_0$ and $\xi_1$ denote the smallest and largest eigenvalues of the
covariance matrix $\Sigma$ such that $1\leq \xi_1/\xi_0 < \infty$ uniformly in
$n$. Then
\begin{equation}\label{eq:eigen_bound}
\xi_0 I_n \leq \Sigma \leq \xi_1 I_n.
\end{equation}
Inequality (\ref{eq:eigen_bound}) is similar in spirit to (\ref{eq:Gram}) and
this similarity will be exploited later on.
%appears in the context of wavelet-vaguelette
%decomposition in a linear inverse problem setting. There, under
%certain regularity conditions, one can associate a pair of vaguelette
%frames $({\cal U}, {\cal V})$ associated to the linear operator $K$ appearing
%in the inverse problem specified by the equation $g = Kf$, where $g$
%is the observed function and $f$ is the true function to be estimated
%from the noisy data. Equations like (\ref{eq:eigen_bound}) describe the
%frame property of the biorthogonal systems ${\cal U}$ and ${\cal V}$
%(w.r.t. certain wavelet basis of sufficient regularity) where the
%matrix $\Sigma$ can be thought of as the Gram matrix of the basis
%elements of the system ${\cal U}$. We shall return to this connection
%later when we shall deal with inverse problems specifically.
% when the parameter vector $\theta$ belongs to a Besov
% sequence space.
% %satisfies an $l^p$ bound of the form $\parallel \theta \parallel_p \leq C$ for
% %some $C > 0$ and we would like to find
% We propose a scheme that achieves the asymptotic minimax rate,  as $\epsilon
% \rightarrow 0$, over such a parameter space when one uses the $l_2$ loss
% function for evaluating the performance of the estimator.

% \subsection{Penalized regression}\label{subsec:penalized_model_selection}

% We first propose the penalized regression estimator that will be applied
% level-wise to the vaguelette coefficients with appropriate choice of the
% parameters. This criterion takes the form

Our primary aim here is to develop a rate-adaptive estimation scheme for the
model (\ref{eq:general}). We use a penalized least squares criterion
\begin{equation}\label{eq:crit}
{\cal R}(\theta,y;\epsilon) = ~\parallel y - \theta \parallel^2 + \epsilon^2
\mbox{pen}(N(\theta)),
\end{equation}
where $\parallel \cdot \parallel$ denotes the $\ell_2$ norm, $N(\theta)$ is the
number of nonzero coordinates of $\theta$, and pen$(\cdot)$ is a nonnegative
function defined on the nonnegative integers.
%Motivated by the optimality properties of \textit{False Discovery Rate (FDR)}
%control procedure of \citeN{AbramovichBDJ06}, we choose our penalty function to
%be
%\begin{equation}\label{eq:penalty}
%\mbox{pen}(k) =  \xi_1 \zeta k (1+\sqrt{2L_{n,k}})^2, \qquad \mbox{where}
%\qquad L_{n,k} = (2\beta + 1) \log(\gamma_n n/k),
%\end{equation}
%with $\zeta > 1$ and $\gamma_n \geq 1$.
%The parameter $q$ can be thought of, in analogy to the FDR procedure as the maximal false discovery rate.
%Also, observe that the standard hard thresholding procedure can be expressed as
%minimizing (\ref{eq:crit}) where $p(k) = c_n k$ for some constant $c_n > 0$.
%Indeed, the value $c_n = \sqrt{2\log n}$ is widely used in sparse signal
%estimation literature (see, for example, \citeN{DonohoJKP97}).
We consider a class of penalty functions of the form
\begin{equation}\label{eq:penalty_general}
\mbox{pen}(k) = \xi_1 \zeta k (1+\sqrt{2L_{n,k}})^2
\end{equation}
where $\zeta > 1$ and $L_{n,k}$ is of the form
\begin{equation}\label{eq:L_n_k}
L_{n,k} = (1+2\beta)\log(\nu n/k), ~~k=1,2,\ldots,n,
\end{equation}
for some $\beta \geq 0$ and $\nu >  e^{1/(1+2\beta)}$ (this condition will be
made clear in Section \ref{sec:upper_bound}). Here $\beta$ is an auxiliary parameter that can be taken to
be zero in the direct estimation problem, but will be positive for the WVD
setting. For now, we treat $\beta$ as a generic parameter taking only
nonnegative values. The choice of the penalty function is motivated by an equivalent
formulation to the False Discovery Rate (FDR) control procedure studied in
\cite{AbramovichBDJ06}. Specifically, for the direct estimation problem (i.e.,
when $\beta=0$), the choice $\nu = 2/w$ with $w \in (0,1)$
corresponds to controlling the FDR at $w$. Qualitatively similar penalties in
the context of direct estimation also appear in \cite{FosterS97}.
\cite{BirgeM01} carried out a systematic study of complexity penalized model
selection in the direct estimation problem, and obtained non-asymptotic bounds
using a penalty class similar to but more general than that used here.

% \subsection{Estimator and an empirical complexity}

% The penalty term pen$(N(\theta))$ in (\ref{eq:crit}) corresponds to a penalty
% on the complexity of the model that we choose to represent $\theta$.
We define our complexity penalized estimator as
\begin{equation}\label{eq:estimator}
\widehat \theta = \arg \min_{\theta} {\cal R}(\theta, y; \epsilon).
\end{equation}
The estimator $\hat \theta$ is given by hard thresholding with a data
dependent threshold.
Indeed, define $\lambda_{n,k} =
\sqrt{\xi_1\zeta}(1+\sqrt{2L_{n,k}})$,
so that the penalty function pen$(k) = k\lambda_{n,k}^2$.
Let $|y|_{(k)}$ denote the order statistics of $|y_i|$:
% Also, let the ordered $y_i$'s (in
% absolute value) be denoted by
$|y|_{(1)} \geq |y|_{(2)} \geq \cdots \geq
|y|_{(n)}$, and let
%  Since $k \mapsto \mbox{pen}(k)$ is a monotone increasing function,
% it can be easily verified that the estimator $\widehat \theta$ is a hard
% thresholding estimator where the threshold is chosen at $\epsilon
% \lambda_{n,\widehat k}$ satisfying the relationship
\begin{equation}\label{eq:threshold}
\widehat k = \arg \min_{k \geq 0} \sum_{i > k} y_{(i)}^2 + \epsilon^2 k
\lambda_{n,k}^2
\end{equation}
Finally, let $t_k^2 = k \lambda_k^2 - (k-1) \lambda_{k-1}^2
= \text{pen}(k) - \text{pen}(k-1)$.
Then it can be shown that $\hat \theta$ is given by hard thresholding
at $t_{\hat k}$ and, for the choices (\ref{eq:penalty_general}) and
(\ref{eq:L_n_k}), that $t_k \approx \lambda_k$ in the sense that
$|t_k - \lambda_k| \leq c/\lambda_k$
\citet[Proposition 11.2 and Lemma 11.7]{Johnstone2013}.

% where, for $k=0$ one can take the sum to be 0 and the corresponding threshold
% to be $\infty$.
%For a fixed $\theta$ one can similarly define theoretical
%complexity $K(\theta',\theta)$ as a function of $\theta'$. Let $\theta_0$ be
%the minimizer of $K(\theta',\theta)$ with respect to $\theta'$.

%Our strategy to bound the risk would be to utilize the basic relationship
%between the theoretical and the empirical complexity
%\begin{equation}\label{eq:basic_bound}
%K(\widehat \theta, \theta) \leq K(\theta_0,\theta) + 2\epsilon \langle z,
%\widehat \theta -\theta_0\rangle .
%\end{equation}

\section{Besov sequence space and the minimax bounds}\label{sec:besov_minimax}

We consider the idealized setting where the parameter $\theta =
(\theta_{jk}:k=1,\ldots,2^j; j=0,1,\ldots) \in \mathbb{R}^\infty$ belongs to a
Besov sequence space determined by a smoothness parameter $\alpha$ and a norm
index $p$. For $\alpha > 1/p - 1/2$, we define the Besov sequence space
$\Theta_{p,q}^\alpha(C)$ for $C > 0$ as
\begin{equation}\label{eq:besov_space}
\Theta_{p,q}^\alpha(C) = \{\theta \in \mathbb{R}^\infty : \sum_{j=0}^\infty
2^{(\alpha - 1/p + 1/2)qj}\parallel \theta_j \parallel_p^q \leq C^q\}
\end{equation}
where $\theta_j = (\theta_{jk})_{k=1}^{2^j}$ and $\parallel \cdot
\parallel_p$ denotes the $\ell_p$ norm.

%\textcolor{red}{[TO DO: Everywhere, change the Besov scale parameter to $q$
%rather than $s$.]}

%\subsection{Penalized estimation under the WVD paradigm}

We estimate $\theta_j = (\theta_{jk} : k=1,\ldots,2^j)$ for $j \geq j_0$ by
applying the penalization method described in Section
\ref{sec:correlated_noise}  separately for each dyadic level $j \geq
j_0$
%\leq j \leq j_1$ where $j_1 = \frac{1}{2\beta+1}\log_2 \epsilon^{-2}$ and
where $j_0$ is an arbitrary but fixed index $\geq 1$. We estimate the
coefficients $\theta_{j_0 k}$ by their empirical value $y_{j_0 k}$. For $j <
j_0$, we set $\widehat\theta_j = y_j$. For $j \geq j_0$, we obtain the
penalized estimator of $\theta_j$ as
\begin{equation}\label{eq:WVD_estimator_level_j}
\widehat \theta_j = \arg\min_{\mu \in \mathbb{R}^{2^j}} \parallel y_j - \mu
\parallel^2 + \epsilon_j^2  \mbox{pen}_j(N(\mu))
\end{equation}
where $\mbox{pen}_j(\mu) = \xi_j \zeta N(\mu) (1 + \sqrt{2L_{n_j,N(\mu)}})^2$
where $N(\mu) =$ number of nonzero coordinates in $\mu$ and $L_{n_j,k}$ is as
in (\ref{eq:L_n_k}), with $n_j$ replacing $n$, and $\nu_{n,j}$ satisfying
\begin{equation}\label{eq:gamma_n_j}
\nu_{n,j} = \begin{cases} \nu  & ~\mbox{if}~ j \leq j_\epsilon \\
\nu [1+(j - j_\epsilon)]^2 & ~\mbox{if}~  j > j_\epsilon \\
\end{cases} \qquad\mbox{where} ~~ j_\epsilon := \log_2 \epsilon^{-2} ~~\mbox{and}~~ \nu > e^{1/(1+2\beta)}.
\end{equation}
As will be shown later, this choice of $\nu_{n,j}$ ensures sufficient control on
the MSE corresponding to each Besov shell, which ensures rate adaptivity of the
proposed estimator.
%We set $\widehat\theta_j = 0$ for $j > j_1$.
For simplicity of exposition, we take $j_0 = 1$ for the rest of the paper since
it does not affect the asymptotic bounds.  From now on we refer to the vector
$(\widehat\theta_{j_0},\widehat\theta_{j_0+1},\ldots)$ as $\widehat \theta$.

Now we state the main contribution of this paper. The most important novelty of
the proposed estimator is that it is rate adaptive, i.e., its rate of
convergence under the squared error loss attains the minimax bound up to a
constant factor over a wide class of Besov sequence spaces. We also show that
the minimax risk for estimation of $\theta$ under the squared error loss
undergoes a phase transition depending on the value of the hyper-parameter
$\boldsymbol{\gamma} = (\alpha,p,q,\beta)$, where $\alpha$, $p$ and $q$
describe the Besov sequence space and the parameter $\beta$ describes the decay
of singular values of the operator. Specifically, as the noise level $\epsilon
\to 0$, there exist three different rate exponents depending on
$\boldsymbol{\gamma}$.
Let $\Gamma_0 = \{ \boldsymbol{\gamma} : \alpha > (1/p - 1/2)_+\}$: this
ensures that $\Theta_{p,q}^\alpha (C)$ is compact in $\ell_2$.
\begin{itemize}
\item[(i)] ``Dense'' regime: \qquad $\Gamma_d :=  \{\boldsymbol{\gamma} : \alpha
> (2\beta+1)\left(1/p - 1/2\right)_+, ~\mbox{and}~p > 0\}$;

\item[(ii)] ``Sparse'' regime: \qquad $\Gamma_s :=
  \{\boldsymbol{\gamma} : \alpha < (2\beta+1)\left(1/p - 1/2\right)~
  \mbox{and}~ 0 < p < 2\} \cap \Gamma_0$;

\item[(iii)] ``Critical'' regime:  \qquad $\Gamma_c := \{\boldsymbol{\gamma} :
\alpha = (2\beta+1)\left(1/p -1/2\right)~ \mbox{and}~ 0 < p < 2\}$,
\end{itemize}
%Our primary aim is to show that in these three different situations the
%proposed penalized least squares estimator has the optimal rate of convergence.
When $0 < p < 2$, we prove the adaptivity of the proposed estimator under an
additional assumption, namely, $\alpha + \beta > 1/p$ (see also Remark \ref{rem:alpha_beta_zone}).
Note that the condition
$\alpha > 1/p$ is necessary for the Besov function class $B_{p,q}^\alpha$ to
embed in spaces of continuous functions (cf. \cite{Johnstone2013}). By the
frame property of vaguelette systems, there exist $0 < \Xi_0 \leq 1 \leq \Xi_1$
such that $\Xi_0 \leq \xi_j \leq \Xi_1$ for all $j \geq 1$. The quantities
$\Xi_0$ and $\Xi_1$ also enter in the minimax bounds even though their roles
are not made explicit. Figures 1 and 2
%\ref{fig:regions_gamma_beta_small} and \ref{fig:regions_gamma_beta_large}
depict the different regimes in the
$(1/p,\alpha)$ plane, for $\beta \in [0,1/2]$ and  $\beta > 1/2$, respectively.
\begin{figure}[th]
\begin{center}
\includegraphics[width=3.6in, height=3.1in]{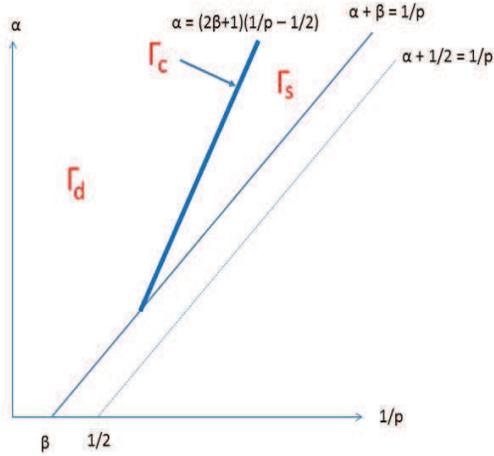}
\caption{Different regimes ($\Gamma_d$: ``dense''; $\Gamma_s$: ``sparse'';
$\Gamma_c$: ``critical'') for the rate of convergence when $0 \leq \beta \leq 1/2$.}
\end{center}
\label{fig:regions_gamma_beta_small}
\end{figure}

\begin{figure}[th]
\begin{center}
\includegraphics[width=3.6in, height=3.1in]{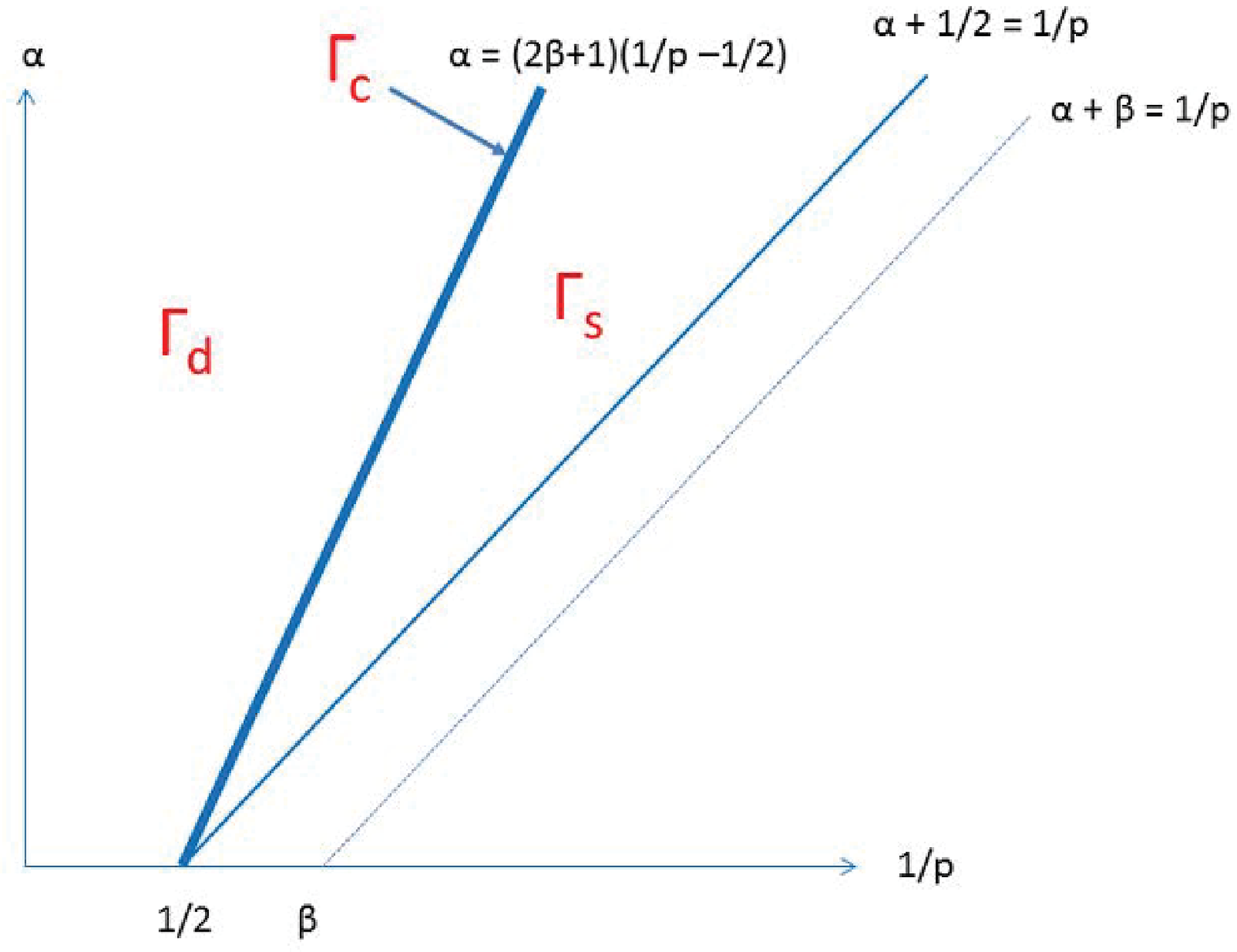}
\caption{Different regimes ($\Gamma_d$: ``dense''; $\Gamma_s$: ``sparse'';
$\Gamma_c$: ``critical'') for the rate of convergence when $\beta > 1/2$.}
\end{center}
\label{fig:regions_gamma_beta_large}
\end{figure}

\begin{theorem}\label{thm:minimax_risk}
Assume observation model (\ref{eq:WVD_transform}) with
$\epsilon_j = \epsilon 2^{\beta j}, \beta \geq 0$.
%Let $\alpha' = \alpha + (1/2-1/p) - (1/2-1/p)_+$ and
Suppose that $\alpha > (1/p - 1/2)_+$ for all $p$,
and $\alpha + \beta > 1/p$ for $0 < p <
2$. Define the rate exponent
\begin{equation}\label{eq:minimax_phases}
r = \begin{cases} \frac{2\alpha}{2\alpha + 2\beta + 1} & \boldsymbol{\gamma} \in \Gamma_d\\
            \frac{2\alpha - 2/p + 1}{2\alpha + 2\beta -2/p + 1} &
            \boldsymbol{\gamma} \in \Gamma_s \\
            1-p/2 & \boldsymbol{\gamma} \in \Gamma_c.
\end{cases}
\end{equation}
Then, for $0 < \epsilon <\epsilon_0$ ($\epsilon_0$ may depend on $C$), we have
\begin{equation}\label{eq:MSE_upper}
\inf_{\widehat \theta} \sup_{\theta \in \Theta_{p,q}^\alpha(C)} \mathbb{E}
\parallel \widehat \theta - \theta\parallel^2  ~\asymp~
\begin{cases}C^{2(1-r)}\epsilon^{2r} & \boldsymbol{\gamma} \in \Gamma_d \\
  C^{2(1-r)}\epsilon^{2r}(1+ \log (C/\epsilon))^r & \boldsymbol{\gamma} \in
  \Gamma_s \\
  C^{2(1-r)}\epsilon^{2r} (1+ \log (C/\epsilon))^{r+(1-p/q)_+} &
  \boldsymbol{\gamma} \in \Gamma_c
\end{cases}
\end{equation}
where ``$\asymp$'' means that both sides are within positive constant multiples
of each other, where the constants depend on $\alpha$, $\beta$,
$p$, $\zeta$, $\Xi_0$ and $\Xi_1$. Moreover, the optimal rates are
attained by the estimator defined through (\ref{eq:WVD_estimator_level_j}) and
(\ref{eq:gamma_n_j}).
\end{theorem}

\begin{remark}\label{rem:alpha_beta_zone}
We can expand the domain of applicability of Theorem \ref{thm:minimax_risk} to
$\{\boldsymbol\gamma : \alpha + \beta > 1/p - 1/2 + 1/(2K)\}$ where $K > 1$, when $0 < p < 2$,
if we also have $\alpha > 1/p - 1/2$. However, this
requires modifying the expression of $j_\epsilon$ in (\ref{eq:gamma_n_j}) to $K \log_2 \epsilon^{-2}$.
\end{remark}

%\textcolor{red}{[TO DO: Re-state Theorem \ref{thm:minimax_risk} by making it
%clear that the optimal rates are attained by the proposed estimator.]}

%\textcolor{red}{[TO DO: Add a note saying that we can push the envelope for
%asymptotic optimality up to the range $\alpha + \beta > 1/p - 1/2$ when $0 < p
%< 2$.]}

\section{Upper bound on the risk}\label{sec:upper_bound}

This section outlines the approach to establishing the upper bounds
in Theorem \ref{thm:minimax_risk}.

% In this section, our goal is to establish the upper bounds in the three cases
% described in Theorem \ref{thm:minimax_risk}. Our approach heavily relies on the
% tools developed in Chapters 11 and 12 of \cite{Johnstone2013}.
% The proof of the rate
% upper bound in the critical regime is given in the Supplementary Material.

%\subsection{Upper bounds for estimation at a single resolution level}
%\label{sec:upper-bounds-estim}

\textbf{Oracle inequalities at a single resolution level}. \
As a first step towards deriving upper bounds on the risk of $\widehat\theta =
(\widehat \theta_j)_{j \geq j_0}$, with $\widehat\theta_j$ defined by
(\ref{eq:WVD_estimator_level_j}), we bound the
risk of the estimator defined by (\ref{eq:estimator}) through the
``oracle inequalities'' that bound the maximal
empirical complexity in terms of the maximal theoretical complexity plus an asymptotically small
term.

With a slight abuse of notation, we write $L_J$ for $L_{n,n_J}$, where $J
\subset \{1,\ldots,n\}$ and $n_J = |J|$. Then define
\begin{equation}\label{eq:M_prime_complexity}
M_n' = \sum_{J \neq \{\}} e^{-L_J n_J},
\end{equation}
where the sum is taken over all subsets $J$ of $\{1,\ldots,n\}$.
% We have the
% following bound (proved in the Supplementary Material).
As long as $\nu_n > e^{1/(1+2\beta)}$, we have,
\begin{equation}
  \label{eq:M_prime_bound}
  M_n' \leq C_{\beta}n^{-2\beta}\nu_n^{-1},
\end{equation}
for some $C_\beta > 0$, as is shown in SM.
% \begin{lemma}\label{lemma:M_prime_bound}
% As long as $\nu_n > e^{1/(1+2\beta)}$, we have,  $M_n' \leq
% C_{\beta}n^{-2\beta}\nu_n^{-1}$ for some $C_\beta > 0$.
% \end{lemma}
% Lemma \ref{lemma:M_prime_bound} will be used to provide a non-asymptotic bound
% on the risk of the penalized regression estimator.

For any $\theta \in \mathbb{R}^n$, let $\theta_J = P_J \theta$, where
$(P_J y)_i =  y_i$ if $i \in J$, and $(P_J y)_i  = 0$, if $i \not\in J$.
%Note that, $P_J$ is the orthogonal projection matrix onto the space $S_J :=
%\mbox{span}\{e_i : i \in J\}$, where $e_i$ is the unit vector in $\mathbb{R}^n$
%with $i$-th coordinate nonzero and the rest zeros.
Now, let $\widehat\theta_J = P_J y$. Then, define the complexity criterion
\begin{equation}\label{eq:compexity_criterion_C}
C_{\epsilon}(J,y) = \parallel y - \widehat \theta_J \parallel^2 + \epsilon^2
\mbox{pen}(n_J) = \sum_{i \not\in J} y_i^2 + \epsilon^2
\mbox{pen}(n_J).
\end{equation}
Define
\begin{equation}\label{eq:optimal_J}
\widehat J = \arg\min_{J \subset \{1,\ldots,n\}} C_{\epsilon}(J,y)
\end{equation}
and observe that,
%based on our discussion above,
$\widehat\theta = P_{\widehat J} y = \widehat\theta_{\widehat J}$. Moreover, if
we define
\begin{equation}\label{eq:K_0}
{\cal R}(\theta,\epsilon) :=
%{\cal R}(\theta_0,\theta;\epsilon) =
\inf_{\theta'}
\parallel \theta - \theta'
\parallel^2 + \epsilon^2 \mbox{pen}(N(\theta')),
\end{equation}
then $\min_J C_{\epsilon}(J,\theta) = {\cal R}(\theta,\epsilon)$.

%In case of $\beta = 0$ one gets {\it Proposition 1} which applies for all
%possible values of $\theta \in \mathbb{R}^n$. However, for future applications,
%we need a slightly more refined result and that requires splitting the domain
%of $\theta$ into two parts, one where the norm is ``small'' and other when the
%norm is ``large''.

%\begin{proposition}\label{prop:oracle_global}
%Let $\widehat \theta$ be the penalized least squares estimator of
%(\ref{eq:crit}) and (\ref{eq:estimator}) for the penalty
%(\ref{eq:penalty_general}) and constant $M$ defined in (\ref{eq:M_complexity}).
%Then there exists a constant $D = D(\zeta)$ such that
%\begin{equation}\label{eq:oracle_global}
%\mathbb{E}\parallel \widehat\theta - \theta \parallel^2 \leq D \left[2M \xi_1
%\epsilon^2 + \min_J C(J,\theta)\right].
%\end{equation}
%The constant $D$ may be taken as $D(\zeta) = 2\zeta(\zeta+1)^3(\zeta-1)^{-3}$.
%\end{proposition}

The next step is the following non-asymptotic bound on the risk of the
penalized least squares estimator (\ref{eq:crit}) which is especially useful
for dealing with our problem. This a restatement of Theorem 11.9
of \cite{Johnstone2013}.
\begin{proposition}\label{prop:oracle_inverse}
Let $\widehat \theta$ be the penalized least squares estimator of
(\ref{eq:crit}) and (\ref{eq:estimator}) for the penalty
(\ref{eq:penalty_general}) and with $M_n'$ defined by
(\ref{eq:M_prime_complexity}). Then there exists a constant $D = D(\zeta)$ such
that
\begin{equation}\label{eq:oracle_inverse}
\mathbb{E}\parallel \widehat\theta - \theta \parallel^2 \leq D\left[2M_n' \xi_1
\epsilon^2 + \min_J C_{\epsilon}(J,\theta)\right] = D\left[2M_n' \xi_1
\epsilon^2 + {\cal R}(\theta,\epsilon)\right].
\end{equation}
The constant $D$ may be taken to be $2\zeta(\zeta+1)^3(\zeta-1)^{-3}$.
\end{proposition}

%\bigskip
We will need to bound the `ideal risk' $\mathcal{R}(\theta, \epsilon)$
over certain $\ell_p$ balls
$\ell_{n,p}(C) = \{x \in \mathbb{R}^n : \sum_{i=1}^n |x_i|^p \leq
C^p\}$.
To state the bound, we introduce \textit{control functions}
$r_{n,p}(C)$.
For $C > 0$ and $0 < p < 2$, let
\begin{equation}\label{eq:r_n_p_C_sparse}
r_{n,p}(C) =
\begin{cases}
C^2 & ~\mbox{if}~~ C \leq \sqrt{1+\log n},\\
C^p[1+\log(n/C^p)]^{1-p/2} & ~\mbox{if}~~ \sqrt{1+\log n} \leq C \leq n^{1/p},\\
n & ~\mbox{if}~~ C \geq n^{1/p},\\
\end{cases}
\end{equation}
while for $p\geq2$, let
\begin{equation}\label{eq:r_n_p_C_dense}
r_{n,p}(C) =
\begin{cases}
n^{1-2/p} C^2 & ~\mbox{if}~~ C \leq n^{1/p},\\
n & ~\mbox{if}~~ C \geq n^{1/p}.\\
\end{cases}
\end{equation}
When $p < 2$, we shall refer to the region $C \geq n^{1/p}$ as the ``dense
zone'', the region $\sqrt{1+\log n} \leq C \leq n^{1/p}$ as the ``sparse zone''
and the region $C \leq \sqrt{1+\log n}$ as the ``highly sparse zone''. When $p
\geq 2$, we shall refer to the region $C \geq n^{1/p}$ as the ``large signal
zone'' and the region $C \leq n^{1/p}$ as the ``small signal zone''.

The proof of the next bound is given in SM.
\begin{lemma}\label{lem:ideal-to-control}
For the `ideal risk' defined in (\ref{eq:crit})-~(\ref{eq:L_n_k}),
there exists $c > 0$ such that
\begin{equation}\label{eq:ideal-to-control}
  \sup_{\theta \in \ell_{n,p}(C)} \mathcal{R}(\theta, \epsilon)
   \leq c (\log \nu) \epsilon^2 r_{n,p}(C/\epsilon).
\end{equation}
\end{lemma}

%\subsection{A general MSE bound}\label{subsec:general_MSE}

\textbf{A general MSE bound.} \
Now we establish a general purpose upper bound for the risk of the estimator
$\widehat \theta$ when $\theta \in \Theta_{p,q}^\alpha(C)$. Let
\begin{equation*}
T(\theta,\epsilon) = \mathbb{E}_\theta\parallel \widehat \theta -  \theta
\parallel^2 = \sum_{j\geq j_0} \mathbb{E}_\theta \parallel \widehat
\theta_j - \theta_j \parallel^2.
\end{equation*}
%By standard bias-variance decomposition, the mean squared error of $\widehat
%\theta$ can be bounded as
%\begin{equation}\label{eq:MSE}
%\mathbb{E}_\theta\parallel \widehat \theta -  \theta \parallel^2 \leq
%T(\theta,\epsilon) + \Delta(\theta),
%\end{equation}
%where $\Delta(\theta)$ is the ``tail bias'' term:
%\begin{equation}
%\Delta(\theta) = \sum_{j=j_1+1}^\infty \parallel \theta_j\parallel^2,
%\end{equation}
%and $T(\theta,\epsilon)$ is the ``variance'' term:
%\begin{equation}\label{eq:variance}
%T(\theta,\epsilon) = \sum_{j=j_0}^{j_1}\mathbb{E}_\theta \parallel \widehat
%\theta_j - \theta_j \parallel^2.
%\end{equation}
By Proposition \ref{prop:oracle_inverse} we have the following bound:
\begin{equation}\label{eq:T_bound}
T(\theta,\epsilon)/D \leq  2\sum_{j \geq j_0} \xi_j M_j' \epsilon_j^2 +
\sum_{j\geq j_0} {\cal R}_{j}(\theta_j,\epsilon_j)
=: T_1(\epsilon) + T_2(\theta, \epsilon),
\end{equation}
say, where $M_j'$ is the analog of $M_n'$ (defined in
(\ref{eq:M_prime_complexity}))
when $n$ is replaced by $n_j$, $\xi_1$ by $\xi_j$, $\nu$ by $\nu_{n_j}$, and
$$
{\cal R}_j(\theta_j,\epsilon_j) := \min_{\theta_j'}
\parallel \theta_j' -\theta_j\parallel^2 + \epsilon_j^2
\mbox{pen}(N(\theta_j'))
$$
is the theoretical complexity in level $j$, and $D > 0$ is some
constant.
The bound (\ref{eq:M_prime_bound}) is constructed to offset the
geometric growth of $\epsilon_j^2 = 2^{2 \beta j} \epsilon^2$ and
together with the choice (\ref{eq:gamma_n_j}) of $\nu_{n_j}$, we
obtain the bound
% Denote the first term within parentheses in (\ref{eq:T_bound}) by
% $T_1(\epsilon)$ and the second term by $T_2(\theta,\epsilon)$, i.e.,
% $T_2(\theta,\epsilon) := \sum_{j \geq j_0} {\cal R}_j(\theta_j,\epsilon_j)$.
%Using (\ref{eq:M_prime_bound}) and (\ref{eq:gamma_n_j}), we have the bound:
\begin{equation}\label{eq:T_1_bound}
T_1(\epsilon) \leq c(\zeta,\Xi_1,\beta,\nu) \epsilon^2 \log \epsilon^{-2},
\end{equation}
which shows that this term is asymptotically negligible.

In order to deal with $T_2(\theta, \epsilon)$, first observe that with $a = \alpha + 1/2 - 1/p$,
\begin{equation}\label{eq:besov_shell}
\theta \in \Theta_{p,q}^\alpha(C) \quad \Longrightarrow \quad
\| \theta_j \|_p \leq C_j := C 2^{-aj},
%\parallel \theta_j \parallel_p \leq C 2^{-(\alpha - 1/p + 1/2)j},
\qquad \forall ~j \geq 1.
\end{equation}
We bound $T_2(\theta, \epsilon)$ by using bounds for the
theoretical complexities $\mathcal{R}_j(\theta_j, \epsilon_j)$ over
the corresponding Besov shells.
Indeed, with $R_j := \epsilon_j^2 r_{n_j,p}(C_j/\epsilon_j)$, from
(\ref{eq:ideal-to-control}) we have
\begin{equation}\label{eq:T_2_bound_preliminary}
\sup_{\theta \in \Theta_{p,q}^\alpha(C)} T_2(\theta,\epsilon)
 =  \sup_{\theta
\in \Theta_{p,q}^\alpha(C)}\sum_{j\geq j_0} {\cal R}_j(\theta_j,\epsilon_j)
\leq c \sum_{j\geq j_0} (\log \nu_{n,j})  R_j.
\end{equation}

\textit{``Dense'' regime:}  \
Here, $\alpha > (2\beta + 1)(1/p - 1/2)_+$ and so $r = 2\alpha/(2\alpha + 2\beta + 1)$.
%Define
%\begin{equation}\label{eq:risk_bound_Besov_shell}
%$R_j = r_{n_j,p}(C_j,\epsilon_j)$.
%\end{equation}
We show that for $p \geq 2$, there exists an index $j_*$---which we
allow to be real valued---such that $R_j$
reaches its peak $R_* = R_{j_*}$ at $j=j_*$ and decays geometrically away from
it. Specifically, we show that
\begin{equation}\label{eq:R_star_dense}
R_* = C^{2(1-r)}\epsilon^{2r}.
\end{equation}
The index $j_*$ is determined by solving the equation
$C_{j_*} = \epsilon_{j_*} n_{j_*}^{1/p}$,
% (without loss of generality, taking $j_*$ to
% be a real number).
% \begin{equation}\label{eq:j_star}
% C_{j_*} = \epsilon_{j_*} n_{j_*}^{1/p},
% \end{equation}
i.e.,   at the ``large signal -- small signal'' boundary (see
(\ref{eq:r_n_p_C_dense})). Note that
%(\ref{eq:j_star})
this equation reduces to
\begin{equation}\label{eq:j_star_expand}
2^{(\alpha + \beta+1/2)j_*} = (C/\epsilon).
\end{equation}
We also show that, for $p \geq 2$,
\begin{equation}\label{eq:Besov_shell_risk_p_geq_2_dense}
R_j = \begin{cases}
R_* 2^{(2\beta + 1)(j-j_*)} &~\mbox{if}~j \leq j_*\\
R_* 2^{-2\alpha(j-j_*)} & ~\mbox{if}~ j \geq j_*.\\
\end{cases}
\end{equation}

For $0 < p < 2$, we have an additional index $j_+ > j_*$ which is obtained from the
equation
$C_{j_+} = \epsilon_{j_+} (1+\log n_{j_+})^{1/2}$,
% \begin{equation}\label{eq:j_plus}
% C_{j_+} = \epsilon_{j_+} (1+\log n_{j_+})^{1/2},
% \end{equation}
i.e., at the ``sparse -- highly sparse zone'' boundary (see
(\ref{eq:r_n_p_C_sparse})).  Thus, $j_+$ satisfies
\begin{equation}\label{eq:j_plus_expand}
2^{(\alpha + \beta -1/p + 1/2)j_+} (1+\log n_{j_+})^{1/2} = (C/\epsilon).
\end{equation}
In this case, there is a second peak of $R_j$ at $j = j_+$. Defining $R_+ =
R_{j_+}$, from (\ref{eq:r_n_p_C_sparse}), we deduce that $R_+ = C^2 2^{-2a
j_+}$. We also show that when $p < 2$,
\begin{equation}\label{eq:Besov_shell_risk_p_lessthan_2_dense}
R_j = \begin{cases}
R_* 2^{(2\beta + 1)(j - j_*)} & ~\mbox{if}~ j < j_* \\
R_* 2^{-p\rho(j-j_*)}[1+\varphi(j-j_*)]^{1-p/2}& ~\mbox{if}~ j_* \leq j < j_+ \\
R_+ 2^{-2a(j-j_+)} & ~\mbox{if}~ j \geq j_+, \\
\end{cases}
\end{equation}
where $\rho := \alpha - (2\beta + 1) (1/p - 1/2) > 0$ and $\varphi = p(\alpha +
\beta + 1/2) \log 2$.
The schematic behavior of the shell risk is depicted in Figure 3.
%\ref{fig:shell_risk_plot}.
In particular, using (\ref{eq:j_star_expand}), (\ref{eq:j_plus_expand})
and (\ref{eq:r_prime_r_sparse_comparison}) (stated below),
it can be checked that $R_* \geq R_+$ for small enough $\epsilon$.
The proofs of (\ref{eq:R_star_dense}) and
(\ref{eq:Besov_shell_risk_p_lessthan_2_dense}) are as in Section 12.5 of
\cite{Johnstone2013}, and hence are given in SM.
From (\ref{eq:j_plus_expand}), we deduce that
\begin{equation}\label{eq:j+disp}
j_+ = \delta^{-1} \log_2(C/\epsilon)(1+o(1))~~~\mbox{as}~\epsilon \to 0,~~~\mbox{where}~~
\delta := \alpha + \beta - 1/p + 1/2.
\end{equation}
Since $\alpha + \beta > 1/p$, so that $\delta > 1/2$, we have $j_+ < j_\epsilon$ for $0 < \epsilon <
\epsilon_0(C)$ (compare with Remark \ref{rem:alpha_beta_zone}). Thus, by the geometric decay of $R_j$ for $j \geq j_+$,  and (\ref{eq:gamma_n_j}) and (\ref{eq:T_2_bound_preliminary}), the risk upper bound follows in the setting $0 < p < 2$. When $p \geq 2$, by (\ref{eq:j_star_expand}) we have $j_* < j_\epsilon$, and so a similar argument, now
involving (\ref{eq:Besov_shell_risk_p_geq_2_dense}), establishes
the risk upper bound.

\textit{``Sparse'' regime:} \
Now, we consider the setting where $0 < p < 2$, $\alpha < (2\beta + 1)(1/p -
1/2)$ and $\alpha + \beta
> 1/p$. The basic strategy is similar to that in the dense case, namely,
bounding $R_j$ by splitting the scale indices $j$'s into three parts: $j \leq
j_*$, $j_* < j < j_+$ and $j \geq j_+$, respectively.
%{\textcolor{red} [in Supplementary Material?]}
%  we can easily
% deduce that $j_+$ satisfies the following bound:
% \begin{equation}\label{eq:j_plus_bounds}
% \frac{1}{\delta} \log_2\left(\frac{C}{\epsilon}\right) - \frac{1}{2\delta}
% \log_2\left(1+\frac{\log 2}{\delta}\log_2\left(\frac{C}{\epsilon}\right)\right)
% < j_+ < \frac{1}{\delta} \log_2\left(\frac{C}{\epsilon}\right) -
% \frac{1}{2\delta} \log_2\left(1+\frac{\log
% 2}{2\delta}\log_2\left(\frac{C}{\epsilon}\right)\right)
% \end{equation}
% provided $1 + (1/\delta)\log 2 \log_2(C/\epsilon) < C/\epsilon$, which happens
% if $\epsilon < \epsilon_0(C)$, say.
Since $R_{j_+} = C^2 2^{-2aj_+}$,
noticing that $r = a/\delta$, by (\ref{eq:j_star_expand}), we have
\begin{equation}\label{eq:R_j_plus_sparse_case}
R_+ = R_{j_+} = C^2 2^{-2aj_+} = C^2 \left(\frac{C^2}{\epsilon^2}\right)^{-r}
(1+\log n_{j_+})^{r} \asymp C^{2(1-r)} \epsilon^{2r}(1+ \log(C/\epsilon))^{r}
\end{equation}
as $\epsilon \to 0$, where the last step follows from (\ref{eq:j+disp}).
%the previous remark about $j_+$.
%(\ref{eq:j_plus_bounds}).

For $j \not\in [j_*,j_+)$, the equalities in
(\ref{eq:Besov_shell_risk_p_lessthan_2_dense}) remain valid, while
it is shown in SM that
\begin{equation}\label{eq:R_j_middle_bound}
R_j \leq R_+ 2^{-\tau(j_+ -j)} ~~\mbox{for}~ j_* \leq j < j_+,
\end{equation}
% We show that,
% \begin{equation}\label{eq:Besov_shell_risk_sparse}
% R_j \leq \begin{cases}
% R_* 2^{(2\beta + 1)(j - j_*)} & ~\mbox{if}~ j < j_* \\
% R_+ 2^{-\tau(j_+ -j)}& ~\mbox{if}~ j_* \leq j < j_+ \\
% R_+ 2^{-2a(j-j_+)} & ~\mbox{if}~ j \geq j_+, \\
% \end{cases}
% \end{equation}
where $\tau = (2\beta + 1) - p(\alpha + \beta + 1/2) = -p[\alpha - (2\beta +
1)(1/p - 1/2)] = - p\rho > 0$.
% The first and the last inequalities in
% (\ref{eq:Besov_shell_risk_sparse}) are in fact equalities and they are
% restatements of the first and last equalities in
% (\ref{eq:Besov_shell_risk_p_lessthan_2_dense}). So, we only need to verify the
% inequality for $j_* < j < j_+ $. This is done in the Supplementary Material.

%\textcolor{red}{[TO DO: Observe that $\tau = -  \rho p$ where $\rho$ is defined
%right after (\ref{eq:Besov_shell_risk_p_lessthan_2_dense}). Need to see if the
%middle bounds in (\ref{eq:Besov_shell_risk_p_lessthan_2_dense}) and
%(\ref{eq:Besov_shell_risk_sparse}) can be proved in a unified way.]}

%define $r' = \alpha/(\alpha + \beta + 1/2)$ and
Observe that
\begin{equation}\label{eq:r_prime_r_sparse_comparison}
\frac{\alpha}{\alpha + \beta + 1/2} \geq \frac{\alpha -1/p + 1/2}{\alpha +
\beta - 1/p + 1/2}  ~~\Leftrightarrow~~ \alpha \leq (2\beta + 1)
\left(\frac{1}{p} - \frac{1}{2}\right),
\end{equation}
while equality on one side implies equality on the other. Defining $r' = \alpha/(\alpha + \beta + 1/2)$,
by (\ref{eq:j_star_expand}),
\begin{equation}\label{eq:R_j_star_sparse_case}
R_* = R_{j_*} = n_{j_*} \epsilon_{j_*}^2 = \epsilon^2 2^{(2\beta+1) j_*} =
\epsilon^2 \left(\frac{C^2}{\epsilon^2}\right)^{1-r'} = C^2 \left(\frac{C^2}{\epsilon^2}\right)^{-r'}.
\end{equation}
Thus, recalling (\ref{eq:R_j_plus_sparse_case}), from
(\ref{eq:r_prime_r_sparse_comparison}) we conclude that $R_{j_*} \leq R_{j_+}$.
Combining, we obtain the result.
Again, since $j_+ < j_\epsilon$ for $0 < \epsilon <
\epsilon_0(C)$, by (\ref{eq:gamma_n_j}) and (\ref{eq:T_2_bound_preliminary}),
the risk upper bound follows.

The proof of the rate upper bound in the \textit{``critical'' regime} is given in SM.

%\textcolor{red}{[TO DO: Add the figures depicting the graphs of $R_j$ vs. $j$
%in the ``dense'' and ``sparse'' regimes.]}

\begin{figure}[th]
\begin{center}
\includegraphics[width=3.8in, height=3in]{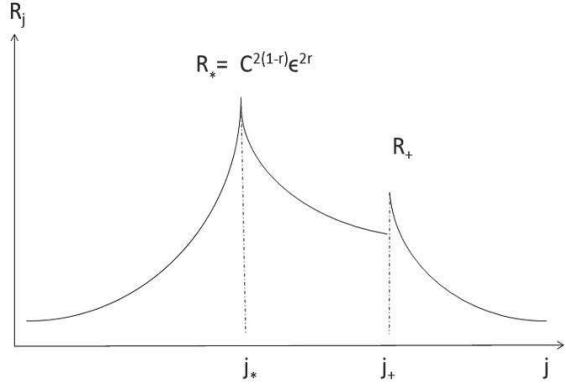}
\caption{Schematic behavior of ``shell risk'' $R_j$, with $j$ treated as a real
variable.}
\end{center}
\label{fig:shell_risk_plot}
\end{figure}

%\subsection{Negligibility of tail bias}\label{subsec:negligibility_tail_bias}

%Recalling (\ref{eq:tail_bias}),
%$$
%\sup_{\theta \in \Theta_{p,s}^\alpha(C)} \Delta(\theta) \leq c C^2
%  (\epsilon^2)^{\frac{2\alpha'}{2\beta+1}}
%\ll C^{2(1-r)}\epsilon^{2r} (1+\log (C/\epsilon))^{r\mathbf{1}\{p<2\}}.
%$$
%We check this carefully in the different regimes. First, consider the ``dense''
%case, i.e., $\alpha
%> (2\beta+1)(1/p - 1/2)$. If $p\geq 2$ then $\alpha' = \alpha$ and so
%$$
%\frac{2\alpha'}{2\beta+1} = \frac{2\alpha}{2\beta+1} >
%\frac{2\alpha}{2\alpha + 2\beta + 1} = r.
%$$
%Now let $p < 2$. Then
%$$
%\frac{2\alpha'}{2\beta+1} = \frac{2(\alpha - 1/p + 1/2)}{2\beta+1} >
%\frac{2\alpha}{2\alpha + 2\beta + 1} = r\quad \Leftrightarrow \quad
%2\alpha(\alpha -1/p + 1/2)
%> (2\beta + 1)(1/p - 1/2)
%$$
%which holds since $\alpha \geq 1/p$.

%Now, suppose $p < 2$ and we consider both the ``sparse'' and ``critical''
%cases, i.e., $\alpha \leq (2\beta + 1)(1/p-1/2)$. Since $\alpha \geq 1/p$, in
%both cases we have
%$$
%\frac{2\alpha'}{2\beta+1} = \frac{2(\alpha - 1/p + 1/2)}{2\beta+1} \geq
%\frac{\alpha - 1/p + 1/2}{\alpha - 1/p + \beta + 1/2} = r.
%$$

%\vskip.1in
%{\it DO WE NEED THE CONDITION THAT $\alpha \geq 1/p$? IT SEEMS TO ME
%  THAT IF WE TAKE $j_1 = \log_2\epsilon^{-2}$ INSTEAD OF $\frac{1}{2\beta+1}
%\log_2\epsilon^{-2}$ AND SIMPLY ASSUME THAT $\alpha > (1/p - 1/2)_+$
%THE RESULTS OF THEOREM 1 SHOULD HOLD. (We only need this when $0 < p <2$).}

\section{Lower bound on the risk}\label{sec:lower_bound}

The idea for the risk lower bound is to minorize the minimax risk of
the model (\ref{eq:WVD_transform}) by the minimax risk of a i.i.d. Gaussian
noise model with covaraince matrix $\Xi_0 I$. Then a lower bound on the latter
is obtained by considering a restricted parameter space for $\theta$ such that
all the level-wise components $\theta_j$ are 0 except for certain specific
dyadic levels $j$, and in those levels the vectors $\theta_j$ are restricted to
lie in $\ell_p$ balls of appropriate radii. Thereafter we can use minimax risk
asymptotics for $\ell_p$ balls \cite{Johnstone2013} to show that the lower bound
thus obtained is of the right asymptotic order.

%\subsection{Equivalence to white noise}

\textbf{Equivalence to white noise.} \
We first show that the minimax risk with noise $z_{jk}$ in the model
(\ref{eq:WVD_transform}) can be bounded below by the minimax risk from a white
noise model. Indeed, let $\Xi_0$ and $\Xi_1$ be as in (\ref{eq:Gram}). Then we
define a new model
\begin{equation}\label{eq:uncorrelated_model}
\tilde y_{jk} = \theta_{jk} + \Xi_0 \epsilon_j  w_{jk}, \qquad k=1,\ldots,2^j,
\quad j \geq j_0,
\end{equation}
where $w_{jk}$ are i.i.d. $N(0,1)$. We denote the minimax risk for estimating
$(\theta)$ in $\ell_2$ loss and with scale parameter $\epsilon$ in model
(\ref{eq:WVD_transform}) by $R_z(\Theta_{p,q}^\alpha(C),\epsilon) :=
\inf_{\widehat\theta}\sup_{\theta \in \Theta_{p,q}^\alpha(C)} \parallel
\widehat\theta - \theta \parallel_2^2$, and that in model
(\ref{eq:uncorrelated_model}) by $R_w(\Theta_{p,q}^\alpha(C),\Xi_0 \epsilon)$.
Then, using Lemma 4.28 of \cite{Johnstone2013}, we conclude that
\begin{equation}\label{eq:minimax_risk_comparison}
R_z(\Theta_{p,q}^\alpha(C),\epsilon) \geq R_w(\Theta_{p,q}^\alpha(C),\Xi_0
\epsilon).
\end{equation}
Thus, it suffices to provide lower bounds on the latter quantity that match
with the bounds in Theorem \ref{thm:minimax_risk}. In the next
subsection, we give an outline of the rate lower bound in the ``dense'' and
``sparse'' regimes.
%  The proof for the critical regime is given in the
% Supplementary Material.

%\subsection{Lower bound : ``dense'' and ``sparse'' cases}\label{subsec:lbd_dense_sparse}

%\textcolor{red}{[TO DO: Shorten this subsection by referring to \cite{Donoho95}
%and the relevant subsections of \cite{Johnstone2013}.]}

\textbf{Lower bound : ``dense'' and ``sparse'' regimes.} \
In both the ``dense'' and ``sparse'' regimes, our strategy is to consider
restricted parameter spaces that are Besov-shells
%\begin{equation}
$\Theta^{(j)}(C) := \{\theta : \parallel\theta_j
\parallel_p \leq C_j ~\mbox{and}~\theta_{j'} = 0~\mbox{if}~j' \neq j\}$,
%\end{equation}
for appropriately chosen $j$. Then, $\Theta^{(j)}(C)$ is isomorphic to the
$\ell_p$ ball $\ell_{n_j,p}(C_j)$. Let $R_w(\Theta^{(j)},\Xi_0\epsilon)$ denote
the minimax risk over $\Theta^{(j)}$, and let $R_N(\ell_{n_j,p}(C_j),\Xi_0
\epsilon_j)$ denote the minimax risk (for estimating $\theta_j$) over the
parameter space $\ell_{n_j,p}(C_j)$, both with respect to the $\ell_2$-loss. Since
$\Theta^{(j)}(C) \subset \Theta_{p,q}^\alpha(C)$, and the $\ell_2$ loss is
coordinate-wise additive, we have
\begin{equation}\label{eq:minimax_risk_lbd_ell_p}
R_w(\Theta_{p,q}^\alpha(C),\Xi_0 \epsilon) \geq
R_w(\Theta^{(j)}(C),\Xi_0\epsilon) \geq R_N(\ell_{n_j,p}(C_j),\Xi_0
\epsilon_j).
\end{equation}

We treat the ``dense'' regime
%i.e. when $\alpha > (2\beta+1)(1/p-1/2)_+$
first. Consider the Besov shell $\Theta^{(j_*)}(C)$, where $j_*$ is
defined by (\ref{eq:j_star_expand})
%(\ref{eq:j_star})
(treating $j_*$ as an integer, for simplicity). Then it
follows from Theorem 13.16 of \cite{Johnstone2013} (restated as Theorem
\ref{thm:minimax_risk_l_p} in SM) that
\begin{equation}\label{eq:minimax_risk_lower_l_p_dense}
R_N(\ell_{n_j,p}(C_{j_*}),\Xi_0 \epsilon_{j_*}) \geq c \Xi_0^2 \epsilon_{j_*}^2
n_{j_*}
\end{equation}
for some $c > 0$, for small enough $\epsilon$. Invoking
(\ref{eq:j_star_expand}), we conclude from (\ref{eq:minimax_risk_comparison}),
(\ref{eq:minimax_risk_lbd_ell_p}) and (\ref{eq:minimax_risk_lower_l_p_dense}) that, for some $c' > 0$,
\begin{equation*}
R_z(\Theta_{p,q}^\alpha(C),\epsilon) \geq c' C^{2(1-r)} \epsilon^{2r},
\end{equation*}
where $r=2\alpha/(2\alpha + 2\beta + 1)$.

In the ``sparse'' regime,
%i.e. when $0 < p < 2$ and $\alpha < (2\beta+1)(1/p-1/2)$,
we consider the Besov-shell $\Theta^{(j_+)}(C)$ where
$j_+$ is defined in (\ref{eq:j_plus_expand}).
%(\ref{eq:j_plus}).
Then, by part (b) of Theorem 13.16 of
\cite{Johnstone2013}, we obtain that
\begin{equation}\label{eq:minimax_risk_lower_l_p_sparse}
R_N(\ell_{n_j,p}(C_{j_+}),\Xi_0 \epsilon_{j_+}) \geq c \Xi_0^2 \epsilon_{j_+}^2
\log n_{j_+}
\end{equation}
for some $c > 0$ and for small enough $\epsilon$. Using
(\ref{eq:j_plus_expand}) and (\ref{eq:j+disp}),
%(\ref{eq:j_plus_bounds}),
from (\ref{eq:minimax_risk_comparison}), (\ref{eq:minimax_risk_lbd_ell_p})
and (\ref{eq:minimax_risk_lower_l_p_sparse}), we conclude that
for some $c' > 0$,
\begin{equation}
R_z(\Theta_{p,q}^\alpha(C),\epsilon) \geq c' C^{2(1-r)} \epsilon^{2r}
(1+\log(C/\epsilon))^r,
\end{equation}
where $r=(2\alpha-2/p+1)/(2\alpha + 2\beta -2/p + 1)$.

Proof of the lower bound in the ``critical'' regime is given in SM.

%\vskip.1in\textcolor{red}{[TO DO: Add acknowledgements.]}

\subsection*{Acknowledgement}

The authors thank Laurent Cavalier for helpful discussions whose untimely death
is deeply regretted. Johnstone's research
is partially supported by NSF grant DMS 0906812, Paul's research is partially
supported by NSF grants DMR 1035468 and DMS 1106690.

% In using BibteX, use wb_stat.bst
%\bibliographystyle{wb_stat}
%\bibliography{bibfile}

\bibliographystyle{wb_stat}
\bibliography{wvd}

\begin{thebibliography}{23}
\newcommand{\enquote}[1]{`#1'}
\providecommand{\natexlab}[1]{#1}
\expandafter\ifx\csname urlstyle\endcsname\relax
  \providecommand{\doi}[1]{doi:\discretionary{}{}{}#1}\else
  \providecommand{\doi}{doi:\discretionary{}{}{}\begingroup
  \urlstyle{rm}\Url}\fi

\bibitem[{Abramovich et~al.(2006)Abramovich, Benjamini, Donoho \&
  Johnstone}]{AbramovichBDJ06}
Abramovich, F, Benjamini, Y, Donoho, D \& Johnstone, IM (2006),
  \enquote{Adapting to unknown sparsity by controlling the false discovery
  rate,} \emph{Annals of Statistics}, \textbf{34}, pp. 584--653.

\bibitem[{Abramovich \& Silverman(1998)}]{AbramovichS98}
Abramovich, F \& Silverman, B (1998), \enquote{Wavelet decomposition approaches
  to statistical inverse problems,} \emph{Biometrika}, \textbf{85}, pp.
  115--129.

\bibitem[{Birg\'{e} \& Massart(2001)}]{BirgeM01}
Birg\'{e}, L \& Massart, P (2001), \enquote{Gaussian model selection,}
  \emph{Journal of European Mathematical Society}, \textbf{3}, pp. 203--268.

\bibitem[{Bissantz et~al.(2007)Bissantz, Hohage, Munk \&
  Ruymgaart}]{BissantzHMR07}
Bissantz, N, Hohage, T, Munk, A \& Ruymgaart, F (2007), \enquote{Convergence
  rates of general regularization methods for statistical inverse problems and
  applications,} \emph{SIAM Journal of Numerical Analysis}, \textbf{45}, pp.
  2610--2636.

\bibitem[{Cai(2002)}]{Cai02}
Cai, TT (2002), \enquote{On adaptive wavelet estimation of a derivative and
  other related linear inverse problems,} \emph{Journal of Statistical Planning
  and Inference}, \textbf{108}, pp. 329--349.

\bibitem[{Cavalier(2008)}]{Cavalier08}
Cavalier, L (2008), \enquote{Nonparametric statistical inverse problems,}
  \emph{Inverse Problems}, \textbf{24}, p. 034004.

\bibitem[{Cavalier \& Golubev(2006)}]{CavalierG06}
Cavalier, L \& Golubev, GK (2006), \enquote{Risk hull method and regularization
  by projections of ill-posed inverse problems,} \emph{Annals of Statistics},
  \textbf{34}, pp. 1653--1677.

\bibitem[{Cavalier et~al.(2004)Cavalier, Golubev, Lepskii \&
  Tsybakov}]{CavalierGLT04}
Cavalier, L, Golubev, GK, Lepskii, O \& Tsybakov, AB (2004), \enquote{Block
  thresholding and sharp adaptive estimation in severly ill-posed inverse
  problems,} \emph{Theory of Probability and its Applications}, \textbf{48},
  pp. 426--446.

\bibitem[{Cavalier et~al.(2002)Cavalier, Golubev, Picard \&
  Tsybakov}]{CavalierGPT02}
Cavalier, L, Golubev, GK, Picard, D \& Tsybakov, AB (2002), \enquote{Oracle
  inequalities in inverse problems,} \emph{Annals of Statistics}, \textbf{30},
  pp. 843--874.

\bibitem[{Cavalier \& Raimondo(2007)}]{CavalierR07}
Cavalier, L \& Raimondo, M (2007), \enquote{Wavelet deconvolution with noisy
  eigenvalues,} \emph{IEEE Transactions on Signal Processing}, \textbf{55}, pp.
  2414--2424.

\bibitem[{Cavalier \& Tsybakov(2002)}]{CavalierT02}
Cavalier, L \& Tsybakov, AB (2002), \enquote{Sharp adaptation for inverse
  problems with random noise,} \emph{Probability Theory and Related Fields},
  \textbf{123}, pp. 323--354.

\bibitem[{Donoho(1995)}]{Donoho95}
Donoho, DL (1995), \enquote{Nonlinear solution to linear inverse problems by
  wavelet-vaguelette decomposition,} \emph{Applied Computational Harmonic
  Analysis}, \textbf{2}, pp. 102--126.

\bibitem[{Donoho et~al.(1997)Donoho, Johnstone, Kerkyacharian \&
  Picard}]{DonohoJKP97}
Donoho, DL, Johnstone, IM, Kerkyacharian, G \& Picard, D (1997),
  \enquote{Universal near minimaxity of wavelet shrinkage,} in
  \emph{Festschrift for Lucien Le Cam}, Springer-Verlag, pp. 183--218.

\bibitem[{Foster \& Stein(1997)}]{FosterS97}
Foster, D \& Stein, R (1997), \enquote{An information theoretic comparison of
  model selection criteria,} Tech. rep., Department of Statistics, University
  of Pennsylvania.

\bibitem[{Johnstone(1999)}]{Johnstone99}
Johnstone, IM (1999), \enquote{Wavelet shrinkage for correlated data and
  inverse problems : adaptivity results,} \emph{Statistica Sinica}, \textbf{9},
  pp. 51--83.

\bibitem[{Johnstone(2013)}]{Johnstone2013}
Johnstone, IM (2013), \emph{Gaussian {E}stimation : {S}equence and {W}avelet
  {M}odels}, Cambridge University Press, manuscript, available at
  http://www-stat.stanford.edu/$\sim$imj/.

\bibitem[{Johnstone et~al.(2004)Johnstone, Kerkyacharian, Picard \&
  Raimondo}]{JohnstoneKPR04}
Johnstone, IM, Kerkyacharian, G, Picard, D \& Raimondo, M (2004),
  \enquote{Wavelet deconvolution in a periodic setting,} \emph{Journal of the
  Royal Statistical Society, Series B}, \textbf{66}, pp. 1--27.

\bibitem[{Kalifa \& Mallat(2003)}]{KalifaM03}
Kalifa, J \& Mallat, S (2003), \enquote{Thresholding estimators for linear
  inverse problems and deconvolutions,} \emph{Annals of Statistics},
  \textbf{31}, pp. 58--109.

\bibitem[{Kolaczyk(1996)}]{Kolaczyk96}
Kolaczyk, ED (1996), \enquote{A wavelet shrinkage approach to tomographic image
  reconstruction,} \emph{Journal of the American Statistical Association},
  \textbf{91}, pp. 1079--1090.

\bibitem[{Loubes \& Lude\~{n}a(2008)}]{LoubesL08}
Loubes, JM \& Lude\~{n}a, C (2008), \enquote{Adaptive complexity regularization
  for linear inverse problems,} \emph{Electronic Journal of Statistics},
  \textbf{2}, pp. 661--677.

\bibitem[{Loubes \& Rivoirard(2009)}]{LoubesR09}
Loubes, JM \& Rivoirard, V (2009), \enquote{Review of rates of convergence and
  regularity conditions for inverse problems,} \emph{International Journal of
  Tomography and Statistics}, \textbf{15}, pp. 349--373.

\bibitem[{Pensky \& Vidakovic(1997)}]{PenskyV97}
Pensky, M \& Vidakovic, B (1997), \enquote{Adaptive wavelet estimator for
  nonparametric density deconvolution,} \emph{Annals of Statistics},
  \textbf{27}, pp. 2033--2053.

\bibitem[{Rochet(2013)}]{Rochet13}
Rochet, P (2013), \enquote{Adaptive hard-thresholding for linear inverse
  problems,} \emph{ESAIM: Probability and Statistics}, \textbf{17}, pp.
  485--499, \doi{10.1051/ps/2012003}.

\end{thebibliography}

\newpage

\setcounter{section}{0} \renewcommand{\thesection}{S.\arabic{section}}
\setcounter{figure}{0} \renewcommand{\thefigure}{S.\arabic{figure}}
\setcounter{equation}{0} \renewcommand{\theequation}{S\arabic{equation}}
\setcounter{theorem}{0} \renewcommand{\thetheorem}{S\arabic{theorem}}

\section*{Supplementary Material}

\subsubsection*{Proof of equation (\ref{eq:M_prime_bound}):}

Using the Stirling's formula bound $k! > \sqrt{2\pi k} k^k e^{-k}$,
\begin{eqnarray*}
M' \leq \sum_{k=1}^n \frac{n^k}{k!} \left(\frac{k}{n\nu}\right)^{k(1+2\beta)}
&\leq& \sum_{k=1}^\infty \frac{1}{\sqrt{2\pi k}}
\left(\frac{k^{2\beta}}{n^{2\beta}}
\frac{e}{\nu^{1+2\beta}}\right)^k \\
&\leq& \frac{1}{n^{2\beta}\nu} \sum_{k=1}^\infty \frac{k^{2\beta} e}{\sqrt{2\pi
k}} \left(\frac{e}{\nu^{1+2\beta}}\right)^{k-1} \leq
\frac{C_{\beta}}{n^{2\beta}\nu}
\end{eqnarray*}
where, in the last step we used the fact that $\nu > e^{1/(1+2\beta)}$.

\subsubsection*{Proof of Lemma \ref{lem:ideal-to-control}:}

Let $|\theta_{(1)}| \geq \cdots \geq |\theta_{(n)}|$ be a decreasing rearrangement of $\theta$.
Then, it is easy to see that
\begin{equation}\label{eq:risk_basic_bound}
{\cal R}(\theta,\epsilon) \leq \sum_{k=1}^n \theta_{(k)}^2 \wedge \epsilon^2\lambda_k^2
= \epsilon^2 \sum_{k=1}^n (\theta_{(k)}/\epsilon)^2 \wedge \lambda_{n,k}^2,
\end{equation}
where $\lambda_{n,k} = \sqrt{\xi_1 \zeta} (1+\sqrt{2(1+2\beta)\log(\nu n/k)})$.

First, consider the case $p\geq 2$. Setting $k=n$ in (\ref{eq:risk_basic_bound}) and noticing that
$\lambda_{n,n}^2 \leq c \log\nu$ for some $c > 0$, we have
\begin{equation}\label{eq:risk_bound_universal}
{\cal R}(\theta,\epsilon) \leq  c(\log \nu)\epsilon^2 n.
\end{equation}
This bound is valid for all values of $C >0$ and actually for all values of $p > 0$.
Moreover, this bound is dominant in particular in the ``dense zone'':
$C/\epsilon \geq n^{1/p}\sqrt{\log \nu}$, in which case the bound reduces to the form
$c(\log \nu) \epsilon^2 r_{n,p}(C/\epsilon)$. Next, by setting $k=0$ in
(\ref{eq:risk_basic_bound}), we have
\begin{equation*}
{\cal R}(\theta,\epsilon) \leq  n\epsilon^2 (n^{-1}\sum_{k=1}^n |\theta_k/\epsilon|^2)
\leq n\epsilon^2 (n^{-1} \sum_{k=1}^n |\theta_k/\epsilon|^p)^{2/p} \leq
\epsilon^2 n^{1-2/p}(C/\epsilon)^2 = \epsilon^2r_{n,p}(C/\epsilon).
\end{equation*}
Clearly, the latter is  bounded by $c(\log \nu) n^{1-p/2} C^2$ which dominates when
$0 < C/\epsilon < n^{1/p}$.

For $p < 2$, we first notice that since $\theta \in \ell_{n,p}(C)$, it implies that $|\theta_{(k)}|
\leq C k^{-1/p}$ for $k=1,\ldots,n$. Therefore, we obtain for all $k \geq 0$,
\begin{equation*}
\sum_{j > k} \theta_{(j)}^2 \leq C^{2-p}(k+1)^{1-2/p} \sum_{j > k} |\theta_{(j)}|^p
\leq C^2 (k+1)^{1-2/p}.
\end{equation*}
Now, invoking this in (\ref{eq:risk_basic_bound}) and setting $k=0$, we have
\begin{equation*}
{\cal R}(\theta,\epsilon) \leq \epsilon^2 (C/\epsilon)^2
\end{equation*}
which is clearly bounded by $(\log \nu) C^2$,  and the latter is of the form $c(\log\nu) \epsilon^2 r_{n,p}(C/\epsilon)$ in the ``sparse zone'': $C \leq \sqrt{1+\log n}$. For the ``dense zone'':
$C \geq n^{1/p}$, we can use the universal bound (i.e., valid for all $p > 0$) given by
(\ref{eq:risk_bound_universal}) and we observe that it is also of the form
$c(\log \nu) \epsilon^2 r_{n,p}(C/\epsilon)$. Thus, it only remains to prove the bound
(\ref{eq:ideal-to-control}) in the case $0 < p < 2$ and $\sqrt{1+\log n} \leq C \leq n^{1/p}$.
The proof of this follows by using an optimization argument as in Section 11.4 of \cite{Johnstone2013} and is omitted.

\subsubsection*{Proof of equations (\ref{eq:R_star_dense}) and
(\ref{eq:Besov_shell_risk_p_geq_2_dense}) :}

To prove (\ref{eq:R_star_dense}), observe that by (\ref{eq:r_n_p_C_dense}),
\begin{eqnarray}\label{eq:R_j_star_dense}
R_* = R_{j_*} &=& n_{j_*} \epsilon_{j_*}^2
~=~ \epsilon^2 2^{(2\beta + 1)j_*} ~=~ C^2 2^{-2\alpha j_*} \qquad (\mbox{by}~(\ref{eq:j_star_expand}))\nonumber\\
&=& C^2 \left(2^{2(\alpha + \beta + 1/2)j_*}\right)^{-r} = C^{2(1-r)}
\epsilon^{2r}.
\end{eqnarray}
Moreover, we have
\begin{equation}\label{eq:R_plus_R_star_ratio_dense}
\log_2(R_+/R_*) = 2(\alpha j_* - a j_+)
\end{equation}
which follows from the first line of (\ref{eq:R_j_star_dense}) and the
expression for $R_+$. From this and (\ref{eq:j_star_expand}) and
(\ref{eq:j_plus_expand}), and using $\rho > 0$, it can be deduced that $R_+
\leq R_*$, which ensures that the final bound on $\sup_{\theta \in
\Theta_{p,q}^\alpha(C)}T_2(\theta,\epsilon)$ is $O(R_*)$.

For the rest of the proof, we note that $C_j/\epsilon_j$ is a monotonically
decreasing sequence in $j$. We first show that
(\ref{eq:Besov_shell_risk_p_geq_2_dense}) holds when $p \geq 2$. First, if $j
\leq j_*$, then we are in the ``large signal zone'', i.e., $C_j/\epsilon_j \geq
n_j^{1/p}$. Hence, $R_j = n_j \epsilon_j^2 = \epsilon^2 2^{(2\beta + 1)j}$.
Hence, the result holds by the first line of (\ref{eq:R_j_star_dense}). Now, if
$j > j_*$, then we are in the in the ``small signal zone'', i.e.,
$C_j/\epsilon_j < n_j^{1/p}$ so that $R_j = n_j^{1-p/2} C_j^2$, from which the
result follows by (\ref{eq:j_star_expand}) and (\ref{eq:R_j_star_dense}).

\subsubsection*{Proof of equation (\ref{eq:Besov_shell_risk_p_lessthan_2_dense}):}

When $p < 2$, we first note that the first bound (i.e., when $j \leq j_*$), and
its proof are exactly the same as in the case $p \geq 2$. The case $j \geq j_+$
corresponds to ``highly sparse zone'', i.e., $C_j/\epsilon_j \leq (1+\log
n_j)^{1/2}$, and hence we have $R_j = C_j^2 = C^2 2^{-2a j}$, which shows, by
comparing with $R_+$, that the result holds in this case. Finally, we turn to
the setting $j_* \leq j < j_+$, i.e., the ``sparse zone''. In this case, define
$\eta_j = (C_j/\epsilon_j) n_j^{-1/p}$. Then,
\begin{equation}\label{eq:R_j_sparse_zone_expand_dense}
R_j = C_j^p \epsilon_j^{2-p} (1+\log(n_j\epsilon_j^p/C_j^p))^{1-p/2} = n_j
\epsilon_j^2 \eta_j^p (1+\log\eta_j^{-p})^{1-p/2}.
\end{equation}
Thus,
\begin{equation}\label{eq:eta_j_power_p}
\eta_j^{-p} =  2^{p(\alpha + \beta + 1/2) (j-j_*)} \left[2^{(\alpha + \beta +
1/2)j_*}/(C/\epsilon)\right]^p = 2^{p(\alpha + \beta + 1/2) (j-j_*)},
\end{equation}
where the second equality is by (\ref{eq:j_star_expand}). Hence, from
(\ref{eq:R_j_sparse_zone_expand_dense}) and the fact that $R_* = R_{j_*} =
n_{j_*} \epsilon_{j_*}^2$, the result follows.

\subsubsection*{Verification of equation (\ref{eq:R_j_middle_bound}) :}

By (\ref{eq:R_j_sparse_zone_expand_dense}) and (\ref{eq:eta_j_power_p}), and
recalling that $\tau = (2\beta + 1) - p(\alpha + \beta + 1/2)$, we have
\begin{eqnarray*}
R_j
%&=& \epsilon^2 \left(\frac{C}{\epsilon}\right)^p 2^{[(2\beta + 1) -
%p(\alpha+\beta+1/2)] j_+} 2^{-\tau(j_+ - j)}  (1+\log \eta_j^{-p})^{1-p/2} \\
&=& \epsilon^2 \left(\frac{C}{\epsilon}\right)^p 2^{[(2\beta + 1) -
p(\alpha+\beta+1/2)] j} (1+\log \eta_j^{-p})^{1-p/2} \\
&=& \epsilon^2 2^{p[(\alpha+\beta - 1/p+1/2) - \alpha + (2\beta + 1)(1/p -
1/2)] j_+} (1+\log n_{j_+})^{p/2} (1+\log \eta_j^{-p})^{1-p/2} 2^{-\tau(j_+ - j)} \\
&=& \epsilon^2 2^{2\beta j_+} (1+\log n_{j_+})^{p/2} (1+\log
\eta_j^{-p})^{1-p/2} 2^{-\tau(j_+ - j)}\\
&=& \epsilon^2 \left(\frac{C^2}{\epsilon^2}\right)^{1-r} (1+\log n_{j_+})^{r -
(1-p/2)} (1+\log \eta_j^{-p})^{1-p/2}  2^{-\tau(j_+ -
j)}\\
&\leq& C^{2(1-r)} \epsilon^{2r} (1+\log n_{j_+})^{r} 2^{-\tau(j_+ - j)} \\
&=& R_+ 2^{-\tau(j_+ - j)},
\end{eqnarray*}
where the second and fourth equalities are due to (\ref{eq:j_plus_expand}) and
the inequality follows from the fact that $\eta_j$'s are decreasing and that
$n_{j_+} \eta_{j_+}^p = (C_{j_+}/\epsilon_{j_+})^p \geq 1$ (since $C_{j_+} = \epsilon_{j_+} (1+\log n_{j_+})^{1/2}$),
while the last equality is due to (\ref{eq:R_j_plus_sparse_case}).

\subsection*{Proof of upper bound in the ``critical'' regime}

%\textcolor{red}{[OPTION : Shift the derivation for the critical regime to the
%supplementary material.]}

Here, $\alpha = (2\beta + 1)(1/p-1/2)$ and $0 < p < 2$. We again consider three
separate blocks : $j \leq j_*$, $j_* < j < j_+$ and $j \geq j_+$. The treatment
of the first and the last block of indices is the same as in the ``sparse''
case above. So, we focus on the middle block.

The conditions $\theta \in \Theta_{p,q}^\alpha(C)$ and $\alpha = (2\beta + 1)
(1/p - 1/2)$ imply that $\alpha - 1/p +1/2 = 2\beta(1/p -1/2)$ so that
\begin{eqnarray}\label{eq:crit_constr}
\sum_{j=j_*}^{j_+} 2^{2\beta(1/p -1/2)qj}\parallel \theta_j
\parallel_p^q &\leq& C^q
~~~\Rightarrow ~~~ \sum_{j=j_*}^{j_+} 2^{(2\beta q/p) j} \parallel \mu_j
\parallel_p^q ~\leq~ \left(\frac{C}{\epsilon}\right)^q,
\end{eqnarray}
where $\mu_j = (\mu_{jk})_{k=1}^{n_j}$ and $\mu_{jk} :=
\epsilon_j^{-1}\theta_{jk}$. In the following, instead of using the omnibus
bound (\ref{eq:T_2_bound_preliminary}) on $T_2(\theta,\epsilon)$ we use the
more direct bound ${\cal R}_j(\theta_j,\epsilon_j) \leq c(\log \nu_{n,j})
\epsilon_j^2 r_{n_j,p}(\parallel \mu_j\parallel_p)$ (for some constant $c > 0$
independent of the parameters $\theta$, $C$ and $\epsilon$), and then, noticing
that $j_+ < j_\epsilon$ so that $\log \nu_{n,j} = \log \nu$,  maximize the sum
$\sum_{j=j_*}^{j_+} \epsilon_j^2 r_{n_j,p}(\parallel \mu_j\parallel_p)$ subject
to (\ref{eq:crit_constr}).

In view of (\ref{eq:crit_constr}),  since $j \geq j_*$, from the fact that $C_{j_*}
= \epsilon_{j_*} n_{j_*}^{1/p}$ so that $C/\epsilon = 2^{(\alpha+\beta+1/2)j_*}$, we obtain
$$
\parallel \mu_j\parallel_p \leq (C/\epsilon) 2^{-(2\beta/p) j} \leq 2^{(\alpha + \beta + 1/2) j}  2^{-(2\beta/p)
j} = 2^{j/p} = n_j^{1/p}.
$$
Thus, from (\ref{eq:r_n_p_C_sparse}), for $j_* \leq j \leq j_+$,
\begin{eqnarray}\label{eq:control_function_bound_level_j_critical}
r_{n_j,p}(\parallel \mu_j \parallel_p) &\leq& c \max\{\parallel \mu_j
\parallel_p^p (1 + \log (n_j /\parallel \mu_j
\parallel_p^{p}))^{1-p/2}, \parallel \mu_j \parallel_p^2 \}
~\leq~ \parallel  \mu_j \parallel_p^p (1+ \log n_j)^{1-p/2},
\end{eqnarray}
where the second inequality follows by noticing that the bound $\parallel \mu_j
\parallel_p^2$ holds only in the ``highly sparse zone'': $\parallel \mu_j
\parallel_p \leq (1+\log n_j)^{1/2}$.
Thus, we consider a majorizing bound for $\sum_{j=j_*}^{j_+} {\cal
R}_j(\theta_j,\epsilon)$ by maximizing
\begin{equation*}
\epsilon^2  (1+ \log n_{j_+})^{1-p/2} \sum_{j=j_*}^{j_+} 2^{2\beta j} \parallel
\mu_j
\parallel_p^p ~~~\mbox{subject to}~ (\ref{eq:crit_constr}).
\end{equation*}
Set $x_j = 2^{2\beta j} \parallel \mu_j\parallel_p^p$, define $x =
(x_j)_{j=j_*}^{j_+}$, and then the optimization problem reduces to
\begin{equation*}
\mbox{maximize}~~~ \epsilon^2 (1+ \log n_{j_+})^{1-p/2}
\parallel x \parallel_1 \qquad \mbox{subject to}
\qquad \parallel x \parallel_{q/p} \leq (C/\epsilon)^p.
\end{equation*}
The value of this maximum is $\epsilon^2 (1+ \log n_{j_+})^{1-p/2} (j_+ -
j_*)^{(1-p/q)_+} (C/\epsilon)^p$. Now, invoking (\ref{eq:j_star_expand}) and
(\ref{eq:j_plus_expand}), we get $j_+ - j_* <
p[2\beta(2\beta+1)]^{-1}\log_2(C/\epsilon)$, and consequently,
\begin{equation}\label{eq:risk_upper_bound_middle_block_critical}
\sum_{j=j_*}^{j_+} {\cal R}_j(\theta_j,\epsilon) \leq c
\epsilon^2(C/\epsilon)^p (1+ \log(C/\epsilon))^{(1-p/2) + (1-p/q)_+}.
\end{equation}
Notice that since $\alpha = (2\beta + 1)(1/p-1/2)$, we have
$$
r=1-p/2 = \frac{\alpha -1/p + 1/2}{\alpha + \beta -1/p + 1/2} =
\frac{\alpha}{\alpha + \beta + 1} = r',
$$
so that from (\ref{eq:R_j_star_sparse_case}) and
(\ref{eq:R_j_plus_sparse_case}), we have $R_{j_*} \leq R_{j_+}$, and the latter
is dominated by the upper bound in
(\ref{eq:risk_upper_bound_middle_block_critical}). Thus, the upper bound for
$T_2(\theta,\epsilon)$ in the critical case follows by combining with the
bounds on $R_j$ for $j \leq j_*$ and $j \geq j_+$.

\subsection*{Details on equations (\ref{eq:minimax_risk_lower_l_p_dense}) and
(\ref{eq:minimax_risk_lower_l_p_sparse})}

Theorem 13.16 of \cite{Johnstone2013}, stated below, states the asymptotic
behavior of the minimax risk of estimation of $\mu \in \mathbb{R}^n$, under the
data model
\begin{equation}\label{eq:data_independent}
y_k = \mu_k + \epsilon_n z_k, \qquad k=1,\ldots,n,
\end{equation}
where $\epsilon_n > 0$ and the random variables $z_k$ are i.i.d. $N(0,1)$. The
minimax risk is calculated using the squared error loss and over the parameter
space $\ell_{n,p}(C_n)$, with $C_n > 0$, i.e.,
\begin{equation}\label{eq:minimax_risk_l_p}
R_N(\ell_{n,p}(C_n),\epsilon_n) = \inf_{\widehat\mu}\sup_{\mu \in
\ell_{n,p}(C_n)}
\parallel \widehat\mu - \mu \parallel_2^2.
\end{equation}
\cite{Johnstone2013} derived the asymptotic expression for
$R_N(\ell_{n,p}(C_n),\epsilon_n)$, as $\epsilon_n \to 0$,  by first deriving an
expression for the Bayes minimax risk in the univariate (i.e., $n=1$) problem,
under the class of univariate priors
\begin{equation*}
\mathfrak{m}_p(\tau) = \{\pi(d\mu) : \int |\mu|^p \pi(d\mu) \leq \tau^p\},
\end{equation*}
so that, with $y \sim N(\mu,\epsilon^2)$, the Bayes minimax risk with respect
to the class $\mathfrak{m}_p(\tau)$ is given by
\begin{equation*}
\beta_p(\tau,\epsilon) = \inf_{\widehat\mu}\sup_{\pi \in \mathfrak{m}_p(\tau)}
B(\widehat\mu,\pi)
\end{equation*}
where $B(\widehat\mu,\pi)$ denotes the Bayes risk of the estimator
$\widehat\mu$ under the squared error loss, with respect to the prior $\pi$.
Proposition 13.4 of \cite{Johnstone2013} states the properties of
$\beta_p(\tau,\epsilon)$, in particular that it is (1) increasing in $p$; (2)
decreasing in $\epsilon$; (3) strictly increasing, concave and continuous in
$\tau^p$; and (4) $\beta_p(\tau,\epsilon) = \epsilon^2
\beta_p(\tau/\epsilon,1)$ and $\beta_p(a\tau,\epsilon) \leq a^2
\beta_p(\tau,\epsilon)$ for all $a \geq 1$.

Furthermore, if we define $\beta_p(\eta) = \beta_p(\eta,1)$, then Theorem 13.7
of \cite{Johnstone2013} states that, as $\eta \to 0$,
\begin{equation*}
\beta_p(\eta) \sim \begin{cases} \eta^2  & ~\mbox{if}~ 2 \leq p \leq \infty,\\
\eta^p(2\log \eta^{-p})^{1-p/2} & ~\mbox{if}~ 0 < p  < 2.\\
\end{cases}
\end{equation*}
Theorem 13.16 of \cite{Johnstone2013}, which summarizes the asymptotic behavior
of $R_N(C_n,\epsilon_n)$, is stated in terms of the function $\beta_p(\eta)$.

\begin{theorem}\label{thm:minimax_risk_l_p}
(Theorem 13.16 of \cite{Johnstone2013}): Introduce the normalized
signal-to-noise ratios
\begin{equation}
\eta_n = n^{-1/p} (C_n/\epsilon_n) \qquad \delta_n = (2\log n)^{-1/2}
(C_n/\epsilon_n).
\end{equation}
For $2 \leq p \leq \infty$, if $\eta_n \to \eta \in [0,\infty]$, then
\begin{equation}\label{eq:R_N_C_epsilon_l_p_dense}
R_N(C_n,\epsilon_n) \sim n\epsilon_n^2 \beta_p(\eta_n).
\end{equation}
For $0 < p <2$,
\begin{itemize}
\item[(a)] if $\eta_n \to \eta \in [0,\infty]$ and $\delta_n \to \infty$ then
again (\ref{eq:R_N_C_epsilon_l_p_dense}) holds.
\item[(b)] If $\eta_n \to 0$ and $\delta_n \to \delta \in [0,\infty)$, then
\begin{equation}\label{eq:R_N_C_epsilon_l_p_sparse}
R_N(C_n,\epsilon_n) \sim \begin{cases} \lambda_n^2 \epsilon_n^2 ([\delta]^p +
\{\delta^p\}^{2/p}), & ~\mbox{if}~\delta > 0,\\
\lambda_n \epsilon_n^2 \delta_n^2, & ~\mbox{if}~\delta  = 0,\\
\end{cases}
\end{equation}
where $[\cdot]$ and $\{\cdot\}$ denote the integer and fractional parts,
respectively, and $\lambda_n = \sqrt{2\log n}$.
\end{itemize}
\end{theorem}
We apply this result with $n = n_{j}$, $C_n = C_{j}$ and $\epsilon_n = \Xi_0
\epsilon_{j}$, for $j=j_*$ in the ``dense'' case and for $j=j_+$ in the
``sparse'' case. It is easy to verify using (\ref{eq:j_star_expand}),
(\ref{eq:j_plus_expand}) and (\ref{eq:j_plus_expand}) that the conditions for
Theorem \ref{thm:minimax_risk_l_p} are satisfied and some elementary
calculations then lead to (\ref{eq:minimax_risk_lower_l_p_dense}) and
(\ref{eq:minimax_risk_lower_l_p_sparse}).

%Then, in the ``dense'' case, $n=n_{j_*}$, $\eta_n = n_{j_*}^{-1/p}
%(C_{j_*}/\Xi_0 \epsilon_{j_*}) \to \Xi_0^{-1}$ and $\delta_n = (2\log
%n_{j_*})^{-1/2} (C_{j_*}/\Xi_0 \epsilon_{j_*})$.

\subsection*{Proof of lower bound in the ``critical'' regime}%\label{subsec:lbd_critical}

%\textcolor{red}{[TO DO: Shift to the supplementary material.]}

Next, we consider the ``critical'' regime. If $p/q \geq 1$, then the lower bound
on the minimax risk for the ``critical'' regime is a continuation of that of the
``sparse'' regime, since $(2\alpha-2/p+1)/(2\alpha + 2\beta -2/p + 1) = 1-p/2$
when $\alpha = (2\beta+1)(1/p - 1/2)$, with $0 < p < 2$. And so, we can use
exactly the same construction as for the ``sparse'' regime in  Section
\ref{sec:lower_bound} to find the lower bound. However, when $0 < p/q <
1$, the lower bound on the minimax risk in the ``critical'' regime has a
discontinuity from that in the ``sparse'' regime and hence we need a different
construction.

We fix two indices $\underline{j} = \lfloor \rho_1 j_*\rfloor$ and $\bar{j} =
\lceil \rho_2 j_* \rceil$ where $1 < \rho_1 < \rho_2 < 2\beta/(2\beta + 1)$ and
$\lfloor x\rfloor$ and $\lceil x \rceil$ are the floor and ceiling functions
(meaning, respectively, the largest integer $\leq$, and smallest integer
$\geq$, $x$). Then, we consider the parameter space
$$
\Theta_{\rho_1,\rho_2}(C) = \{\theta : \sum_{j=\underline{j}+1}^{\bar{j}} 2^{aq
j} \parallel \theta_j \parallel_p^q \leq C^q~~\mbox{and} ~~\theta_{j} =
0~\mbox{if}~j \leq \underline{j}~\mbox{or}~j > \bar{j}\}.
$$
Clearly, $\Theta_{\rho_1,\rho_2}(C) \subset \Theta_{p,q}^\alpha(C)$ and
therefore,
\begin{equation}\label{eq:minimax_risk_lbd_critical}
R_w(\Theta_{p,q}^\alpha(C),\Xi_0 \epsilon) \geq
R_w(\Theta_{\rho_1,\rho_2}(C),\Xi_0 \epsilon),
\end{equation}
where $R_w(\Theta_{\rho_1,\rho_2}(C),\Xi_0 \epsilon)$ denotes the minimax risk
over $\Theta_{\rho_1,\rho_2}(C)$ under $\ell_2$ loss based on the data from model
(\ref{eq:uncorrelated_model}).

We adopt a Bayes-minimax approach to find a lower bound for
$R_w(\Theta_{\rho_1,\rho_2}(C),\Xi_0 \epsilon)$. Specifically,  following the
construction in Lemma 11 of \cite{DonohoJKP97}, for each $j \in
\{\underline{j}+1,\ldots,\overline{j}\}$, we construct a prior $\Pi_j$ as
follows. For appropriately chosen $n_{0j}$ ($\leq n_j$) and $\delta_{0j} > 0$,
set $\tau_j = n_{0j}/(2n_j)$. Then $\theta_j \sim \Pi_j$ means that the random
variables $\theta_{jk}$, $k=1,\ldots,n_j$, are i.i.d. according to the
distribution which puts mass $1-\tau_j$ at 0 and mass $\tau_j/2$ each at $\pm
\delta_{0j}$. Moreover, we choose the priors $\Pi_j$ to be independent for
different $j \in \{\underline{j}+1,\ldots,\overline{j}\}$. Define restricted
parameter spaces
$$
\Theta_{0j}(n_{0j},\delta_{0j}) = \{\theta_j \in \mathbb{R}^{n_j} :
\#\{\theta_{jk} \neq 0\} \leq n_{0j}, ~~\mbox{and}~~|\theta_{jk}|\leq
\delta_{0j}~\mbox{for all}~k\},
$$
and the restricted priors $\bar{\Pi}_j(\cdot) =
\Pi_j(\cdot|\Theta_{0j}(n_{0j},\delta_{0j}))$ for
$j=\underline{j}+1,\ldots,\bar{j}$. Now, suppose that we can choose
${(n_{0j},\delta_{0j})}_{j=\underline{j}+1}^{\bar{j}}$ in such a way that the
following conditions hold.
\begin{itemize}
\item[(i)]
The set $\{\theta : \theta_j \in \Theta_{0j}(n_{0j},\delta_{0j})~\mbox{for}~j
\in \{\underline{j}+1,\ldots,\bar{j}\}, ~\mbox{and}~\theta_j
=0~\mbox{otherise}\}$ is contained in $\Theta_{\rho_1,\rho_2}(C)$.
\item[(ii)]
There exist $d \in (0,1)$, $d' \in (0,d)$ and an $A > 0$, such that $n_{0j}
\leq A n_j^{(1-d)}$ and $\delta_{0j} \leq \Xi_0 \epsilon_j \sqrt{2(d-d') \log
n_j} $ for all $j=\underline{j}+1,\ldots,\bar{j}$.
\end{itemize}

If (ii) holds, then we proceed as in the proof of Lemma 11 of
\cite{DonohoJKP97}, which uses the bound
\begin{equation}\label{eq:sign_error_bound}
\parallel \widehat \theta_j - \theta_j\parallel_2 \geq
(\delta_{0j}/2) \sum_{k=1}^{n_j} \mathbf{1}(|\widehat\theta_{jk} -\theta_{jk}|
> \delta_{0j}/2),
\end{equation}
derives the form of the univariate Bayes estimator $\widehat\theta_{jk}^*$ for
$\theta_{jk}$ with loss function $\mathbf{1}(|\widehat\theta_{jk} -\theta_{jk}|
> \delta_{0j}/2)$, and then uses large deviations bound for Binomial random
variables to bound the deviation probabilities under $\bar{\Pi}_j$ of the
random variable on the RHS of (\ref{eq:sign_error_bound}) when
$\widehat\theta_j = \widehat\theta_j^*$. From these, we conclude that, there
exists a constant $b
> 0$, not depending on $j$, such that for any estimator $\widehat\theta$ and for each $j \in \{\underline{j}+1,\ldots,\bar{j}\}$,
\begin{equation*}
\mathbb{P}_{\bar{\Pi}_j}(\parallel \widehat\theta_j - \theta_j\parallel_2^2
\geq n_{0j} \delta_{0j}^2 /40) \geq 1-2e^{-b n_{0j}},
\end{equation*}
where $\mathbb{P}_{\bar{\Pi}_j}$ denotes the joint probability of $(\tilde
y,\theta)$ computed under $\bar{\Pi}_j$. Hence, for any $\widehat\theta$,
\begin{equation}\label{eq:risk_prob_lbd_critical}
\mathbb{P}_{\prod_{j=\underline{j}+1}^{\bar{j}}\bar{\Pi}_j}\left(\sum_{j=\underline{j}+1}^{\bar{j}}\parallel
\widehat\theta_j - \theta_j\parallel_2^2 \geq \frac{1}{40}
\sum_{j=\underline{j}+1}^{\bar{j}} n_{0j} \delta_{0j}^2\right) \geq
1-2\sum_{j=\underline{j}+1}^{\bar{j}} e^{-b n_{0j}}.
\end{equation}
Since $\bar{\Pi}_j$ is supported on $\Theta_{0j}(n_{0j},\delta_{0j})$, now
invoking property (i) and using Chebyshev's inequality we conclude from
(\ref{eq:risk_prob_lbd_critical}) that, for small enough $\epsilon$,
\begin{equation}\label{eq:risk_prob_lbd_critical_inf}
\inf_{\widehat\theta} \sup_{\theta \in \Theta_{\rho_1,\rho_2}(C)}
\mathbb{E}\parallel \widehat\theta - \theta
\parallel^2 \geq c \sum_{j=\underline{j}+1}^{\bar{j}} n_{0j} \delta_{0j}^2
\end{equation}
for some $c > 0$, provided
\begin{equation}\label{eq:lbd_prob_tail}
\sum_{j=\underline{j}+1}^{\bar{j}} e^{-b n_{0j}} \to 0 ~~\mbox{as}~~ \epsilon
\to 0.
\end{equation}
We choose $\delta_{0j} = c_0 \Xi_0 \epsilon_j (\log_2(C/\epsilon))^{1/2}$ and
\begin{eqnarray*}
n_{0j} &=& c_1 (C/\epsilon)^p 2^{-2\beta j} (\bar{j}-\underline{j})^{-p/q}
(\log_2(C/\epsilon))^{-p/2} \\
&=& c_1 (\bar{j}-\underline{j})^{-p/q} (\log_2(C/\epsilon))^{-p/2}
2^{j(1-(2\beta+1)(1-j_*/j))},
\end{eqnarray*}
for some constants $c_0, c_1 > 0$. The second expression for $n_{0j}$ follows
from (\ref{eq:j_plus_expand}) and the fact that $\alpha + \beta + 1/2 =
(2\beta+1)/p$. Since $0 < \rho_1 < \rho_2 < (2\beta+1)/2\beta$, it easily
follows that, by choosing $c_0, c_1 > 0$ appropriately, we can ensure that (i),
(ii) and (\ref{eq:lbd_prob_tail}) are satisfied. Finally,
\begin{equation*}
\sum_{j=\underline{j}+1}^{\bar{j}} n_{0j} \delta_{0j}^2 \geq c_2  \epsilon^2
(C/\tilde \epsilon)^p (\log_2(C/\tilde \epsilon))^{(1-p/2)+(1-p/q)},
\end{equation*}
for sufficiently small $\epsilon$, which, together with
(\ref{eq:minimax_risk_comparison}), (\ref{eq:minimax_risk_lbd_critical}) and
(\ref{eq:risk_prob_lbd_critical_inf}) yields the lower bound in Theorem
\ref{thm:minimax_risk} for the ``critical'' regime.

\end{document}